\documentclass{amsart}
\usepackage{amssymb,url}
\usepackage[bookmarks=false,colorlinks=true]{hyperref}
\input{xy}
\xyoption{all}
\usepackage[margin=1.5in]{geometry}

\newtheorem{theorem}{Theorem}[section]
\newtheorem*{theorem*}{Theorem}
\newtheorem*{regtheorem*}{Theorem~\ref{regular-result}}
\newtheorem*{divtheorem*}{Theorem~\ref{divisibility-theorem}}
\newtheorem*{semiregtheorem*}{Theorem~\ref{semiregular-result}}
\newtheorem{lemma}[theorem]{Lemma}
\newtheorem{proposition}[theorem]{Proposition}
\newtheorem{corollary}[theorem]{Corollary}

\theoremstyle{definition}
\newtheorem{definition}[theorem]{Definition}
\newtheorem{example}[theorem]{Example}
\newtheorem{question}[theorem]{Question}
\newtheorem{remark}[theorem]{Remark}

\newcommand{\Z}{\ensuremath{\mathbb{Z}}}
\def\ZZ{{\mathbb Z}}
\def\QQ{{\mathbb Q}}
\def\RR{{\mathbb R}}
\DeclareMathOperator{\coker}{{coker}}
\def\im{{\operatorname{im}}}              
\def\Hom{{\operatorname{Hom}}}           
\def\line{{\operatorname{line \,\,}}}       
\def\sd{{\operatorname{sd \,}}}           
\def\coker{{\operatorname{coker}}}        
\def\lcm{{\operatorname{lcm}}}             
\def\Syl{{\operatorname{Syl}}}             
\def\tetrahedron{{\operatorname{tetra}}}             
\def\octahedron{{\operatorname{octa}}}             
\def\cube{{\operatorname{cube}}}             
\def\dodecahedron{{\operatorname{dodeca}}}             
\def\icosahedron{{\operatorname{icosa}}}             
\def\dcube{{\operatorname{d-cube}}}

\begin{document}

\author[Berget]{Andrew Berget} \email{berget@math.ucdavis.edu}
\author[Manion]{Andrew Manion} \email{andymanion@gmail.com}
\author[Maxwell]{Molly Maxwell} \email{molly.maxwell@coloradocollege.edu}
\author[Potechin]{Aaron Potechin} \email{potechin@mit.edu}
\author[Reiner]{Victor Reiner} \email{reiner@math.umn.edu}
\address{
  Department of Mathematics\\ 
  University of California\\
  Davis, CA 95616
}
\address{
  Department of Mathematics\\
  Princeton University\\
  Princeton, NJ 08540
}
\address{
  Department of Mathematics and Computer Science\\
  Colorado College\\ 
  14 E. Cache la Poudre St.\\
  Colorado Springs, CO 80903
}
\address{
  Department of Mathematics\\
  Massachusetts Institute of Technology\\
  Cambridge, MA 02139
}
\address{
  School of Mathematics\\
  University of Minnesota\\
  206 Church St. S.E.\\
  Minneapolis, MN 55455
}

\thanks{Work of all authors supported by NSF grants DMS-0245379 and
DMS-0601010, and completed partly during REU programs at the University
of Minnesota during the summers of 2003, 2004 and 2008.}

\title{The critical group of a line graph}


\begin{abstract}
  The critical group of a graph is a finite abelian group whose order
  is the number of spanning forests of the graph. This paper provides
  three basic structural results on the critical group of a line
  graph.
  \begin{itemize}
  \item The first deals with connected graphs containing no cut-edge.
    Here the number of independent cycles in the graph, which is known
    to bound the number of generators for the critical group of the
    graph, is shown also to bound the number of generators for the
    critical group of its line graph.
  \item The second gives, for each prime $p$, a constraint on the
    $p$-primary structure of the critical group, based on the largest
    power of $p$ dividing all sums of degrees of two adjacent
    vertices.
  \item The third deals with connected graphs whose line graph is
    regular.  Here known results relating the number of spanning trees
    of the graph and of its line graph are sharpened to exact
    sequences which relate their critical groups.
\end{itemize}
The first two results interact extremely well with the third.  For
example, they imply that in a regular nonbipartite graph, the critical
group of the graph and that of its line graph determine each other
uniquely in a simple fashion.
\end{abstract}

\maketitle

\tableofcontents

\section{Introduction and main results} 
\label{intro}

The critical group $K(G)$ of a graph $G$ is a finite abelian group
whose order is the number $\kappa(G)$ of spanning
forests\footnote{Throughout this paper all spanning forests are
  assumed to be maximal in the sense that adding an edge of $G$ to a
  spanning forest creates a cycle. This removes the possibility of a
  connected graph containing a disconnected spanning forest.}  of the
graph.  One can define $K(G)$ in several ways, closely related to the
{\it cycle} and {\it bond} spaces of the graph, the {\it graph
  Laplacian}, as well as a certain {\it chip-firing game} that is
played on the vertices of the graph and is called the {\it abelian
  sandpile model} in the physics literature.  The interested reader
can find some of the standard results on $K(G)$ in \cite{Bacher,
  Biggs1} and \cite[Chapter 13]{Godsil}.  Some of this material is
reviewed in Sections~\ref{rods-section} and
\ref{graph-critical-groups} below, along with unpublished results from
the bachelor's thesis of D. Treumann \cite{Treumann} on functoriality
for critical groups.

The critical group $K(G)$ and its relation to the structure of the
graph $G$ remain, in general, mysterious.  The goal of this paper is
to compare the structure of the critical group of a {\it connected
  simple graph} (that is, a connected graph having no multiple edges
and no loops) with that of the critical group of its {\it line graph}.
Recall that for a graph $G=(V,E)$, its line graph $\line G=(V_{\line
  G}, E_{\line G})$ has vertex set $V_{\line G}:=E$, the edge set of
$G$, and an edge in $E_{\line G}$ corresponding to each pair of edges
in $E$ that are incident at a vertex.
\begin{figure}[h]
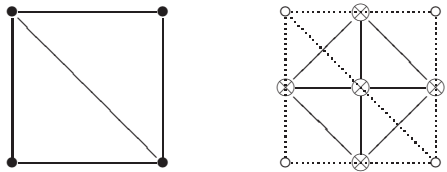

  \[
  \xygraph{
    !{<0cm,0cm>;<1cm,0cm>:<0cm,1cm>::}
    !{(-1,1) }*{\bullet}="u"
    !{(1,-1) }*{\bullet}="v"
    !{(1,1) }*{\bullet}="w"
    !{(-1,-1) }*{\bullet}="x"
    "u"-"w"-"v"-"x"-"u"
    "u"-"v"
  }
  \qquad\qquad 
  \xygraph{
    !{<0cm,0cm>;<1cm,0cm>:<0cm,1cm>::}
    !{(-1,1) }*{\circ}="u"
    !{(1,-1) }*{\circ}="v"
    !{(1,1) }*{\circ}="w"
    !{(-1,-1) }*{\circ}="x"
    !{(0,0) }*{\otimes}="uv"
    !{(0,1) }*{\otimes}="uw"
    !{(0,-1) }*{\otimes}="vx"
    !{(1,0) }*{\otimes}="wv"
    !{(-1,0) }*{\otimes}="xu"
    "u"-@{..}"w"-@{..}"v"-@{..}"x"-@{..}"u"
    "u"-@{..}"v"
    "uw"-"wv"-"vx"-"xu"-"uw"-"uv"
    "uv"-"wv"
    "uv"-"vx"
    "uv"-"xu"
  }
  \]
  \caption{A graph $G$ and its line graph $\line G$ with $G$
    underlayed.}
  \label{fig:linegraph}
\end{figure}
Our main results say that, under three different kinds of hypotheses,
the structure of $K(\line G)$ is not much more complicated than that
of $K(G)$, as we now explain.

\subsection{The hypothesis of no cut-edge}
It is well-known, and follows from one of the definitions of $K(G)$ in
Section~\ref{graph-critical-groups}, that for a connected graph $G$, the number
$\beta(G):=|V|-|E|+1$ of independent cycles in $G$ gives an upper
bound on the number of generators required for $K(G)$; that is,
\begin{equation}
  \label{critical-group-invariant-factor-form}
  K(G) = \bigoplus_{i=1}^{\beta(G)} \ZZ_{d_i}
\end{equation}
where $\ZZ_d$ denotes the cyclic group $\ZZ/d\ZZ$ ({\it not} the
$d$-adic integers), and the $d_i$ are positive integers (some of which
may be $1$).  Our first main result asserts that the same bound on the
number of generators holds for $K(\line G)$ when one assumes that $G$
is {\it $2$-edge-connected}, that is, $G$ is connected and contains no
cut-edge.

%
%

\begin{theorem}
\label{2-edge-connected-theorem}
When the simple graph $G$ is $2$-edge-connected, the critical group
$K(\line G)$ can be generated by $\beta(G)$ elements.
\end{theorem}
\noindent Note that one needs {\it some} hypothesis on the graph $G$
for this conclusion to hold.  For example, a {\it star graph}
$K_{1,n}$ (= one vertex of degree $n$ connected to $n$ vertices of
degree one) has $\beta(K_{1,n})=0$. However, its line graph is the
complete graph $K_n$ and thus, according to
Proposition~\ref{Lorenzini's-theorem} below, has critical group
$K(\line K_{1,n})=\ZZ_n^{n-2}$, requiring $n-2$ generators.
Theorem~\ref{2-edge-connected-theorem} is proven in
Section~\ref{2-edge-connected-section}, using a useful presentation of
$K(\line G)$ given in Section~\ref{line-graph-presentation-section}.

\subsection{The hypothesis that degree sums of adjacent vertices are
  divisible by $p$}
As $K(\line G)$ is a finite abelian group, its structure is completely
determined once one knows, for each prime $p$, the structure of its {\it
  $p$-primary component} or {\it $p$-Sylow subgroup} $\Syl_p(K(\line
G))$.  Section~\ref{divisibility-section} below proves the following
stringent constraint on this $p$-primary structure, based on the
largest power $k(p)$ such that $p^{k(p)}$ divides all of the sums
$\deg_G(v)+\deg_G(w)$ as one runs through all edges $e=\{v,w\}$ in the
edge set $E$ of $G$. Here $\deg_G(v)$ is the number of edges of $G$ with
$v$ as an endpoint; it is the {\it degree of the vertex $v$}.

\begin{theorem}
  \label{divisibility-theorem}
  Let $G=(V,E)$ be a connected simple graph that contains at least one
  cycle of even\footnote{This even length cycle need not be {\it
      minimal}.  For example, a connected graph with two cycles $C_1,
    C_2$ of odd length will also contain a cycle of even length that
    traverses $C_1$, follows a path from $C_1$ to $C_2$, then
    traverses $C_2$ and follows the same path back to $C_1$.} length.
  Use the abbreviated notation $K:=K(\line G)$, and let $p$ be a prime
  for which the quantity $k(p) \geq 1$.

  Then for $G$ bipartite, one has
  $$
  K /p^{k(p)} K 
  \cong \ZZ_{p^{k(p)}}^{\beta(G)-1} \oplus \ZZ_{\gcd(p^k, |V|)},
  $$
  while for $G$ nonbipartite, one has
  $$
  K /p^{k(p)} K
  \cong \ZZ_{p^{k(p)}}^{\beta(G)-2} \oplus 
  \begin{cases}
    0 &\text{ if }p\text{ is odd,} \\
    \ZZ_2^2 &\text{ if }p=2\text{ and }|V|\text{ is even,} \\
    \ZZ_4 &\text{ if }p=2\text{ and }|V|\text{ is odd.} \\
  \end{cases}
  $$
\end{theorem}

\subsection{The regularity hypothesis}
Our third class of main results deals with the situation where $\line
G$ is regular, that is, all its vertices have the same degree.  Say
that a graph is {\it $d$-regular} if all of its vertices have
degree $d$.  It is an easy exercise to check that, for
connected graphs $G$, one has $\line G$ regular only in these two
situations:
\begin{itemize}
\item $G$ itself is $d$-regular.  In this case, $\line G$ will be
  $(2d-2)$-regular.
\item $G$ is bipartite and $(d_1,d_2)${\it -semiregular}, meaning that
  its vertex bipartition $V=V_1 \sqcup V_2$ has all vertices in $V_i$
  of degree $d_i$ for $i=1,2$.  In this case, $\line G$ will be
  $(d_1+d_2-2)$-regular.
\end{itemize}
Two classical theorems of graph spectra  explain how the the numbers of spanning trees $\kappa(G)$ and 
$\kappa(\line G)$ determine each other in this situation. The first is due
originally to Vahovskii \cite{Vahovskii} and later Kelmans \cite{Kelmans}, then rediscovered by Sachs \cite[\S 2.4]{CvetkovicDoobSachs},
while the second is due originally to Cvetkovi\'c \cite[\S
5.2]{Rubey} (see also \cite[Theorem 3.9]{Mohar}). \vskip .1in
\noindent

\begin{theorem*}
Let $G$ be a connected graph with $\line G$ regular.\\
\noindent (\textbf{Sachs}) If $G$ is $d$-regular, then
\begin{equation}\label{Sachs'-theorem}
  \kappa(\line G)= d^{\beta(G)-2} \,\, 2^{\beta(G)} \,\, \kappa(G).
\end{equation}
\noindent (\textbf{Cvetkovi\'c}) If $G$ is bipartite and
$(d_1,d_2)$-semiregular, then
\begin{equation}\label{Cvetkovic's-theorem}
  \kappa(\line G) = 
  \frac{(d_1+d_2)^{\beta(G)}}{d_1 d_2} 
  \left( \frac{d_1}{d_2} \right)^{|V_2|-|V_1|} \,\, \kappa(G). 
\end{equation}
\end{theorem*}
\noindent These results suggest a close relationship between the
critical groups $K(G)$ and $K(\line G)$ in both of these situations.

\subsubsection{Regular graphs}
We focus first on such a relation underlying Sachs' equation
\eqref{Sachs'-theorem}, as here one can be quite precise.

The occurrence of the factor $2^{\beta(G)} \,\, \kappa(G)$ within
\eqref{Sachs'-theorem} suggests consideration of the {\it edge
  subdivision} graph $\sd G$, obtained from $G$ placing a new vertex
at the midpoint of every edge of $G$.  It is well-known that
\begin{equation}
\label{subdivision-complexity-relation}
\kappa(\sd G) = 2^{\beta(G)} \,\, \kappa(G),
\end{equation}
due to an obvious $2^{\beta(G)}$-to-$1$ surjective map from the
spanning trees of $\sd G$ to those of $G$\footnote{More explicitly,
  there are exactly $\beta(G)$ edges that do not lay on a given
  spanning tree of $G$. Upon subdividing, there are $2^{\beta(G)}$ ways
  to extend the resulting tree to a spanning tree of $\sd G$.}.
Underlying this relation, Lorenzini \cite{Lorenzini} observed that the
critical groups $K(\sd G)$ and $K(G)$ also determine each other in a
trivial way: $K(G)$ has the form given in
\eqref{critical-group-invariant-factor-form} if and only if for the
same positive integers $d_1,d_2,\ldots,d_{\beta(G)}$ one has the
following form for $K(\sd G)$:
\begin{equation}
  \label{subdivision-invariant-factor-form}
  K(\sd G) = \bigoplus_{i=1}^{\beta(G)} \ZZ_{2d_i}.
\end{equation}
See Proposition~\ref{Lorenzini's-theorem} below.  In light of
\eqref{subdivision-complexity-relation}, one might expect that
equation \eqref{Sachs'-theorem} generalizes to a short exact sequence
of the form
\begin{equation}
\label{naive-hope}
0 
\rightarrow \ZZ_d^{\beta(G)-2} 
\rightarrow K(\line G)
\rightarrow K(\sd G)
\rightarrow 
0
\end{equation}
where $\ZZ_d$ denotes a cyclic group of order $d$.  This is never far
from the truth. After reviewing and developing some theory of critical
groups and their functoriality in Sections~\ref{rods-section} and
\ref{graph-critical-groups} below, we use functoriality to prove the
following result in Section~\ref{regular-result-section}.
\begin{theorem}
\label{regular-result}
For any connected $d$-regular simple graph $G$ with $d \geq 3$, there
is a natural group homomorphism $ f:K(\line G) \rightarrow K(\sd G) $
whose kernel-cokernel exact sequence takes the form
$$
0 
\rightarrow \ZZ_d^{\beta(G)-2} \oplus C
\rightarrow K(\line G)
\overset{f}{\rightarrow} K(\sd G)
\rightarrow C
\rightarrow 
0
$$
in which the cokernel $C$ is the following cyclic $d$-torsion group:
$$
C=
\begin{cases}
  0        & \text{ if } G \text{ is non-bipartite and }d\text{ is odd},\\
  \ZZ_2 & \text{ if } G \text{ is non-bipartite and }d \text{ is even},\\
  \ZZ_d & \text{ if } G \text{ is bipartite}.
\end{cases}
$$
\end{theorem}

It turns out that Theorems~\ref{2-edge-connected-theorem} and
\ref{divisibility-theorem} interact very well with
Theorem~\ref{regular-result}.  When $G$ is a $d$-regular simple
$2$-edge-connected graph, Theorem~\ref{2-edge-connected-theorem}
implies that $K(\line G)$ needs at most $\beta(G)$ generators, while
Proposition~\ref{Lorenzini's-theorem} implies that $K(\sd G)$ requires at
least $\beta(G)$ generators, forcing $K(\line G)$ to require either
$\beta(G)-1$ or $\beta(G)$ generators.  This shows that the exact
sequence in Theorem~\ref{regular-result} is about as far as possible
from being split, and gives it extra power in determining the
structure of $K(\line G)$ given that of $K(G)$ (and hence also $K(\sd
G)$).

Even more precisely, it will be shown in
Section~\ref{regular-nonbipartite-section} that when $G$ is both
$d$-regular and {\it nonbipartite},
Theorems~\ref{divisibility-theorem} and \ref{regular-result} combined
show that $K(G)$ and $K(\line G)$ determine each other uniquely in the
following fashion.

\begin{corollary}
  \label{regular-nonbipartite-corollary}
  For $G$ a simple, connected, $d$-regular graph with $d \geq 3$ which
  is nonbipartite, after uniquely expressing
  $$
  K(G) \cong \bigoplus_{i=1}^{\beta(G)} \ZZ_{d_i}
  $$ 
  with $d_i$ dividing $d_{i+1}$, one has
  $$
  K(\line G) \cong \left( \bigoplus_{i=1}^{\beta(G)-2} \ZZ_{2dd_i}
  \right) \oplus 
  \begin{cases}
    \ZZ_{2d_{\beta(G)-1}} \oplus \ZZ_{2d_{\beta(G)}} &
    \text{ for }|V|\text{ even,}\\
    \ZZ_{4d_{\beta(G)-1}} \oplus \ZZ_{d_{\beta(G)}} & 
    \text{ for  }|V|\text{ odd.}
  \end{cases}
  $$
\end{corollary}

\subsubsection{Semiregular bipartite graphs}

Section~\ref{semi-regular} uses functoriality to prove the following
result analogous to Theorem~\ref{regular-result} and suggested by
Cvetkovi\'c's equation \eqref{Cvetkovic's-theorem}.
\begin{theorem}
  \label{semiregular-result}
  Let $G$ be a connected bipartite $(d_1,d_2)$-semiregular graph $G$.
  Then there is a group homomorphism 
  $$
  K(\line G) \overset{g}{\rightarrow} K(G)
  $$
  whose kernel-cokernel exact sequence
  \begin{equation}
0 
\rightarrow \ker(g)
\rightarrow K(\line G)
\overset{g}{\rightarrow} K(G)
\rightarrow \coker(g)
\rightarrow 
0
\end{equation}
has $\coker(g)$ all $\lcm(d_1,d_2)$-torsion, and has
$\ker(g)$ all $\frac{d_1+d_2}{\gcd(d_1,d_2)} \lcm(d_1,d_2)$-torsion.
\end{theorem}

\noindent
Note that this result describes the kernels and cokernels less
completely than Theorem~\ref{regular-result}.
Section~\ref{semi-regular} discusses examples illustrating why this is
necessarily the case.

Section~\ref{example-section} illustrates some of the preceding
results by showing how they apply to the examples of complete graphs and
complete bipartite graphs, as well as the $1$-skeleta of
$d$-dimensional cubes and the Platonic solids.

\subsubsection{Directed line graphs}
The reader should compare our results with recent results of Levine
\cite{Levine} on the critical group of a \textit{directed} line
graph. If $G$ is a directed graph, then the directed line graph
$\mathcal{L}G$ is defined so that a pair of directed edges $e$ and $f$
of $G$ are adjacent (and oriented from $e$ to $f$) if the head of $e$
is equal to the tail $f$.

The critical group of a directed graph is defined as the torsion
subgroup of the cokernel of the Laplacian matrix of $G$. Levine proves
\cite[Theorem~1.2]{Levine} that if $G$ is a strongly connceted
Eulerian directed graph, then there is a surjective group homomorphism
from the critical group of $\mathcal{L}G$ to the critical group of
$G$. Moreover, when the $G$ is balanced and $d$-regular,  the
kernel of this homomorphism is the $d$-tosion subgroup of the critical
group of $\mathcal{L}G$. As can be seen from
Theorem~\ref{regular-result}, we do not obtain such easily stated
results in the undirected case.

\section{Some theory of lattices} \label{rods-section}
This section recalls some of the theory of rational lattices in
Euclidean spaces and their determinant groups, along with
functoriality and Pontrjagin duality for these groups, borrowing
heavily from Bacher, de la Harpe, and Nagnibeda \cite{Bacher} and
Treumann \cite{Treumann}.  In the next section, these constructions
will be specialized to critical groups of graphs.

\subsection{Rational orthogonal decompositions}
Consider $\RR^m$ with its usual inner product $\langle \cdot, \cdot
\rangle$ in which the standard basis vectors $e_1,\ldots,e_m$ are
orthonormal.  The $\ZZ$-span of this basis is the integer lattice
$\ZZ^m$.

\begin{definition}
  A {\it rational orthogonal decomposition} is an orthogonal
  $\RR$-vector space decomposition of $\RR^m = B^\RR \oplus Z^\RR$ in
  which $B^\RR$, $Z^\RR $ are $\RR$-subspaces which are rational, that
  is, spanned by elements of $\ZZ^m$.
\end{definition}

\begin{example}\label{graphs}
  The main example of interest for us will be the following, discussed
  further in Section~\ref{graph-critical-groups}.  If $G=(V,E)$ is a
  graph with $|E|=m$, then the space $Z^\RR$ of $1${\it-cycles}
  together with its orthogonal complement, the space $B^\RR$ of {\it
    bonds} or $1${\it-coboundaries}, give a rational orthogonal
  decomposition $\RR^E \cong \RR^m = B^\RR \oplus Z^\RR$.  Here one
  must fix an (arbitrary) orientation of the edges in $E$ in order to
  make the identification $\RR^E \cong \RR^m$.  In the remaining
  sections, the basis element of $\RR^E$ corresponding to an edge
  $\{u,v\}$ of $G$ oriented from $u$ to $v$ will sometimes be denoted
  $e$ and sometimes $(u,v)$, with the convention that
  $(v,u)=-(u,v)=-e$ in $\RR^E$.
\end{example}

An $r$-dimensional rational subspace $\Lambda^\RR \subset \RR^m$
inherits the inner product $\langle \cdot, \cdot \rangle$.  The space
$\Lambda^\RR$ contains two lattices of rank $r$, namely $\Lambda:=
\Lambda^\RR \cap \ZZ^m$ and its {\it dual lattice}
\[
\Lambda^\# := \{ x \in \Lambda^\RR: \langle x, \lambda \rangle \in \ZZ
\text{ for all } \lambda \in \Lambda \}.
\]
Since $\langle \Lambda, \Lambda \rangle \subset \langle \ZZ^m , \ZZ^m
\rangle = \ZZ$, one has an inclusion $\Lambda \subset \Lambda^\#$.
Their quotient is called the {\it determinant group}
$$
\det(\Lambda):= \Lambda^\# / \Lambda.
$$

Given a rational orthogonal decomposition $\RR^m = B^\RR \oplus
Z^\RR$, one obtains two determinant groups $\det(B), \det(Z)$, which
turn out to be both isomorphic to what we will call the {\it critical
  group}
$$
K:=\ZZ^m / (B \oplus Z)
$$
of the rational orthogonal decomposition. Indeed, if $\pi_B, \pi_Z$
denote the orthogonal projections from $\RR^m$ onto $B^\RR, Z^\RR$,
then these maps turn out to give rise to surjections from $\ZZ^m$ onto
$B^\#$ and $Z^\#$, respectively, and which induce isomorphisms (see
\cite[Proposition 3]{Bacher})
$$
\begin{matrix}
  \det(B) & \cong & K & \cong & \det(Z) \\
  B^\#/B & \overset{\pi_B}{\longleftarrow} &\ZZ^m/(B \oplus Z)
  &\overset{\pi_Z}{\rightarrow} & Z^\#/Z .
\end{matrix}
$$ 

One can compute the critical group $K$ very explicitly as the
(integer) cokernel of several matrices, for example via their {\it
  Smith normal form}.  If the lattices $B, Z$ have $\ZZ$-bases $\{
b_1,\ldots,b_\alpha\} , \, \{ z_1, \ldots, z_\beta \}$ then let $M_B, \, M_Z, \,
M_{B \oplus Z}$ be matrices having columns given by
$\{b_i\}_{i=1}^\alpha, \, \{z_j\}_{j=1}^\beta, \, \{b_i\}_{i=1}^\alpha \cup
\{z_j\}_{j=1}^\beta$, respectively.  The {\it Gram matrices} $M_B^t
M_B, \, M_Z^t M_Z$ express the bases for $B, Z$ in terms of the dual
bases for $B^\#, Z^\#$, and hence
\begin{align*}
  K  &  \cong \ZZ^m /(B \oplus Z)   = \coker M_{B \oplus Z}, \\
  & \cong B^\#/B   = \coker (M_B^t M_B),  \\
  & \cong Z^\#/Z = \coker (M_Z^t M_Z).
\end{align*}

\subsection{Functoriality}
\label{functoriality-subsection}
Suppose that one has two rational orthogonal decompositions $\RR^{m_i}
= B^\RR_i \oplus Z^\RR_i$ for $i=1,2$, and an $\RR$-linear map $f:
\RR^{m_1} \rightarrow \RR^{m_2}$.  When does $f$ induce a
homomomorphism $f: K_1 \rightarrow K_2$ between their critical groups?

It is natural to assume that $f$ carries the integer lattice
$\ZZ^{m_1}$ into $\ZZ^{m_2}$, that is, $f$ is represented by a matrix
in $\ZZ^{m_2 \times m_1}$.  Note that this already implies that the
adjoint map $f^t: \RR^{m_2} \rightarrow \RR^{m_1}$ with respect to the
standard inner products will also satisfy $f^t(\ZZ^{m_2}) \subset
\ZZ^{m_1}$, since this map is represented by the transposed $\ZZ^{m_1
  \times m_2}$ matrix.

What one needs further to induce homomorphisms of critical groups is
that $f(B_1) \subset B_2$ and $f(Z_1) \subset Z_2$.  The following
proposition gives a useful reformulation.

\begin{proposition}
\label{Treumann-characterization}
For a linear map $f: \RR^{m_1} \rightarrow \RR^{m_2}$ satisfying
$f(\ZZ^{m_1}) \subset \ZZ^{m_2}$, one has
$$
f(B_1) \subset B_2 \Longleftrightarrow f^t(Z_2) \subset Z_1
\Longleftrightarrow f(Z_1) \subset Z_2^\#
$$
and
$$
f(Z_1) \subset Z_2 \Longleftrightarrow f^t(B_2) \subset B_1
\Longleftrightarrow f(B_1) \subset B_2^\#.
$$
\end{proposition}
\begin{proof}
  All of the implications follow using the adjointness of $f, f^t$
  with respect to the pairings on $\RR^{m_1}, \RR^{m_2}$, along with
  the definitions of $B_i^\#, Z_i^\#$ and the fact that $B_i=
  Z_i^\perp$.
\end{proof}

When a linear map $f: \RR^{m_1} \rightarrow \RR^{m_2}$ satisfies all
of the conditions in the previous proposition, we say that $f$ is a {\it
  morphism of rational orthogonal decompositions}.  It is clear that
$f$ then induces a homomorphism $K_1 \rightarrow K_2$ between the
critical groups, denoted here also by $f$.

Note the following property of such maps $f$ for future use.
\begin{proposition}\label{f-commutes-with-projections}
  Any morphism $f: \RR^{m_1} \rightarrow \RR^{m_2}$ of rational
  orthogonal decompositions intertwines the projection maps onto
  either $B_i^\RR$ or $Z_i^\RR$.  That is, the following diagram
  commutes:
  \begin{equation}
    \xymatrix{
      \RR^{m_1}  \ar[d]_{\pi_{B_1}}\ar[r]^f& \RR^{m_2}\ar[d]^{\pi_{B_2}} \\
      B_1^\RR  \ar[r]_{f^t}& B_2^\RR
    }
  \end{equation}
  and the same holds replacing $B_i$ by $Z_i$ everywhere.
\end{proposition}
\begin{proof}
  Given $x_1 \in \RR^{m_1}$ and $b_2 \in B_2^\RR$, note that
  $$
  \langle \pi_{B_2}(f(x_1)), b_2 \rangle =
  \langle f(x_1) , b_2 \rangle =
  \langle x_1, f^t(b_2) \rangle =
  \langle \pi_{B_1}(x_1), f^t(b_2) \rangle =
  \langle f(\pi_{B_1}(x_1)), b_2 \rangle.
  $$
  Since this equality holds for any test vector $b_2 \in B_2^\RR$, one
  concludes that $\pi_{B_2}(f(x_1)) = f(\pi_{B_1}(x_1))$.
\end{proof}

\subsection{Pontrjagin duality}
Every finite abelian group $K$ is isomorphic to its {\it Pontrjagin
  dual}
$$
K^* := \Hom_\ZZ(K,\QQ/\ZZ).
$$
This isomorphism is not, in general, natural (although the isomorphism
$K \cong K^{**}$ is).  However, for critical groups $K=\ZZ^m/(B \oplus
Z)$ associated with a rational orthogonal decomposition, the
isomorphism comes about naturally from the pairing
$$
\begin{matrix}
\ZZ^m \times \ZZ^m & \rightarrow  & \QQ \\
(x , y )           & \mapsto      & \langle \pi(x),\pi(y) \rangle
\end{matrix}
$$
where $\pi$ is either of the orthogonal projections $\pi_B$ or $\pi_Z$.  
This induces a pairing
$$
\langle \cdot, \cdot \rangle: K \times K \rightarrow  \QQ/\ZZ
$$ 
which is nondegenerate in the sense that the following  map is an isomorphism:
\begin{equation}
\label{Pontrjagin-duality-iso}
\begin{matrix}
K & \rightarrow & \Hom_\ZZ(K,\QQ/\ZZ) &(= K^*)\\
x & \mapsto     & \langle x, \cdot \rangle.
\end{matrix}
\end{equation}

Pontrjagin duality is contravariant in the following sense.  Given a
homomorphism $f: K_1 \rightarrow K_2$ of abelian groups, there is a
dual morphism $f^* : K_2^* \rightarrow K_1^*$ given by $f^*(g) = g
\circ f.$ The next proposition asserts that this duality interacts as
one would expect with morphisms of rational orthogonal decompositions.

\begin{proposition}
\label{Pontrjagin-dual-is-adjoint}
For a morphism $f: \RR^{m_1} \rightarrow \RR^{m_2}$ of rational
orthogonal decompositions, Pontrjagin duality identifies $f^t$ with
$f^*$, in the sense that the following diagram commutes:
\begin{align}
  \label{Pontrjagin-duality-diagram}
  \xymatrix{
    K_2 \ar[d]\ar[r]^{f^t} & K_1\ar[d] \\
    K_2^*\ar[r]_{f^*} & K_1^*
  }
\end{align}
Here the vertical maps are both Pontrjagin duality isomorphisms as in
\eqref{Pontrjagin-duality-iso}.
\end{proposition}
\begin{proof}
  Unravelling the definitions, this amounts to checking that if $x_i
  \in \ZZ^{m_i}$ for $i=1,2$, then one has $\langle f^t(x_2), x_1
  \rangle = \langle x_2, f(x_1) \rangle$.
\end{proof}

This last proposition has a useful consequence.
\begin{corollary}
  \label{identify-kernel-cokernel}
  For a morphism $f: \RR^{m_1} \rightarrow \RR^{m_2}$ of rational
  orthogonal decompositions, the maps induced by $f, f^t$ on critical
  groups satisfy $\ker(f)^* \cong \coker(f^t)$ and $\coker(f)^* \cong
  \ker(f^t)$.
\end{corollary}
\begin{proof}
  Pontrjagin duality generally gives $\ker(f)^* \cong \coker(f^*)$ and
  $\coker(f)^* \cong \ker(f^*)$, so this follows from
  Proposition~\ref{Pontrjagin-dual-is-adjoint}.
\end{proof}

\section{The critical group of a graph} 
\label{graph-critical-groups}

This section particularizes the discussion of critical groups from the
previous section to the context of Example~\ref{graphs}, that is, the
critical group $K(G)$ for a graph $G=(V,E)$.  It also recalls how one
can use spanning trees/forests to be more explicit about some of these
constructions, and reviews for later use some other known results
about critical groups of graphs. 

We will use the term ``spanning tree'' when discussing connected
graphs and ``spanning forest'' when no connectivity is
assumed.

\subsection{Cycles, bonds, Laplacians, and spanning trees}
\label{cycles-bonds-section}
Let $G=(V,E)$ be a graph.  After picking an orientation for its edges,
the usual {\it cellular boundary map} from $1$-chains to $0$-chains
with real or integer coefficients
$$
\begin{aligned}
  \RR^E \overset{\partial_G}{\longrightarrow} \RR^V  \\
  \ZZ^E \overset{\partial_G}{\longrightarrow} \ZZ^V
\end{aligned}
$$
is defined $\RR$- or $\ZZ$-linearly as follows: A basis element $e$
corresponding to an edge directed from vertex $u$ to vertex $v$ is
sent to $\partial_G(e)=+v-u$.  One considers the negative $-e$ of this
basis element as representing the same edge but directed from $v$ to
$u$, which is consistent with
$$
\partial_G(-e) = +u-v = -\partial_G(e).
$$

Elements in the kernel of $Z^\RR := \ker \partial_G$ are called {\it
  cycles}, while elements in the perpendicular space
$B^\RR:=\im \partial_G^t$ are called {\it bonds}.  Thus $\RR^E = B^\RR
\oplus Z^\RR$ is a rational orthogonal decomposition associated with
the graph $G=(V,E)$, and we denote by $K(G)$ the associated critical
group.

The lattice $B$ of bonds is known to be spanned by the signed
incidence vectors $b(V_1,V_2)$ of the directed edges that span across
a {\it cut} (partition) $V=V_1 \sqcup V_2$.  The lattice $Z$ of cycles
is known to be spanned by the signed incidence vectors $z(C)$ coming
from {\it directed cycles} in $G$.

If one wants a smaller $\ZZ$-spanning set for $B$, one can take the
vectors $b_G(\{v\}, V \setminus \{v\})$ for cuts that isolate single
vertices; this vector $b_G(\{v\}, V \setminus \{v\})$ is exactly the
row vector of the $|V| \times |E|$ boundary map $\partial_G$ indexed
by $v$. To simplify notation, we will write
$$
b_G(v):=b_G(\{v\}, V \setminus \{v\}) 
$$
for this bond, and we will call it the bond at $v$ in $G$. In order to
select out of this spanning set a $\ZZ$-basis for $B$, one should omit
exactly one vertex from each connected component of $G$.

Here are a few consequences of these facts:
\begin{enumerate}\label{bond-consequences}
\item[(i)] The Gram matrix $M_B^t M_B$ corresponding to the above
  mentioned $\ZZ$-basis for $B$ gives what is usually called a {\it
    (reduced) Laplacian matrix} $\overline{L(G)}$; the matrix $M_B$ is
  obtained from $\partial_G^t$ by removing the columns corresponding
  to the chosen vertex in each connected component of $G$.  As a
  consequence, one has by Kirchhoff's {\it Matrix-tree Theorem} (see,
  e.g., \cite[Theorem 2.2.12]{West}) that
  $$
  |K(G)| = \det\overline{L(G)} = \kappa(G),
  $$
  the number of spanning forests in $G$.
\item[(ii)]
(The chip-firing/dollar-game/sandpile/Picard presentations for $K(G)$)

Given a connected graph $G=(V,E)$ with boundary map $\ZZ^E
\overset{\partial}{\longrightarrow} \ZZ^V$, bond lattice
$B:=\im \partial^t$, and any vertex $v_0$ in $V$, one has an
isomorphism
$$
\begin{aligned}
  K(G) &\cong \coker (M^t_B M_B)\\
  &\cong \coker (\overline{L(G)})\\
  &\cong \ZZ^{V \setminus \{v_0\}} / \overline{L(G)}\\
  &\cong \ZZ^V / \left(\ZZ v_0 + \im \left(\partial_G \partial_G^t \right) \right)\\
  &\cong \ZZ^V / \left(\ZZ v_0 + \partial_G(B) \right)\\
  &\cong \ZZ^V / \left(\ZZ v_0 + \ZZ(\partial_G \, b_G(v)_{v \in V} )\right).
\end{aligned}
$$

\item[(iii)] For any vertex $u$ of $G$, one has the relation 
  $$ \sum_{\{v\in V: \{u,v\} \in E\}} (u,v)= 0 $$ 
  in $K(G)=\ZZ^E/(B \oplus Z)$.

\item[(iv)] For any directed cycle $u_0 \rightarrow u_1 \rightarrow
  \cdots \rightarrow u_{\ell-1} \rightarrow u_{\ell}=u_0$ in $G$, one
  has the relation $$\sum_{i=0}^{\ell-1} (u_i,u_{i+1})= 0$$ in
  $K(G)=\ZZ^E/(B \oplus Z)$.
\end{enumerate}

Fixing a particular spanning forest $T$ for $G$ allows one to
simultaneously construct $\Z$-bases of $B$ and $Z$. Removing any edge
$e$ in the forest $T$ creates a new connected component in the forest,
say with vertex set $V_e \subset V$; ranging over all edges $e$ in
$T$, the signed incidence vectors $b^T_e$ for the cuts $V=V_e \sqcup(V
- V_e) $ form a $\Z$-basis for $B$. Dually, adding any edge $e$ in $E
- T$ to $T$ creates a unique cycle in $T \cup \{e\}$; ranging over all
edges in $E - T$, the signed incidence vectors $z^T_e$ of these cycles
form a $\ZZ$-basis for $Z$.

\subsection{A presentation for $K( \line
  G)$}\label{line-graph-presentation-section}

Proposition~\ref{alternate-presentation} below gives a useful
presentation for $K(\line G)$ that is an immediate consequence of the
last equation in assertion \eqref{bond-consequences}(ii) above.  It
will be used both in the proof of
Theorem~\ref{2-edge-connected-theorem} and in the analysis of $K(\line
K_n)$ in Section~\ref{complete-graph-section}.

Let $G=(V,E)$ be a connected simple graph, so that $\line
G=(V_{\line G}, E_{\line G})$ is also connected.  Identify the vertex
set $V_{\line G}$ of the line graph of $G$ with the edge set $E$ of
$G$.  After picking arbitrary orientations for the edges of $\line G$,
consider the boundary map for $\line G$:
\[
\partial_{\line G}: \ZZ^{E_{\line G}} \longrightarrow
\ZZ^E \,\, ( = \ZZ^{V_{\line G}}).
\]

\begin{proposition}\label{alternate-presentation}
  Given a connected simple graph $G=(V,E)$ and any edge $e_0$ in $E$,
  one has an isomorphism
  \[
  \begin{aligned}
    K(\line G) & \cong \ZZ^E /(\ZZ e_0 + \partial_{\line G}(B_{\line G}))\\
    &\cong \ZZ^E / \left(\ZZ e_0 + \ZZ(\partial_{\line G} \, b_{\line
        G}(e)_{e \in E}) \right).
\end{aligned}
\]
\end{proposition}

\subsection{Lorenzini's result on edge subdivisions}
\label{subdivision-example}

The {\it edge subdivision} of a graph $G$ is the graph $\sd G$
obtained by creating a new midpoint vertex called $uv$ for every edge
$\{u,v\}$ of $G$; that is, $\{u,v\}$ is removed and replaced by two
edges $\{u, uv\}, \{v,uv\}$ in $\sd G$. An orientation of $G$ induces
an orientation of $\sd G$: If $u$ is oriented towards $v$ then $u$ is
oriented towards $uv$ and $uv$ is oriented towards $v$.

\begin{figure}
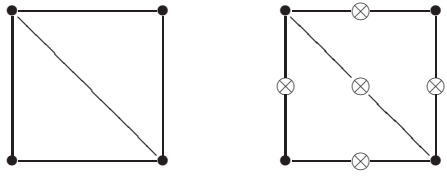

\[
  \xygraph{
    !{<0cm,0cm>;<1cm,0cm>:<0cm,1cm>::}
    !{(-1,1) }*{\bullet}="u"
    !{(1,-1) }*{\bullet}="v"
    !{(1,1) }*{\bullet}="w"
    !{(-1,-1) }*{\bullet}="x"
    "u"-"w"-"v"-"x"-"u"
    "u"-"v"
  }
  \qquad\qquad 
  \xygraph{
    !{<0cm,0cm>;<1cm,0cm>:<0cm,1cm>::}
    !{(-1,1) }*{\bullet}="u"
    !{(1,-1) }*{\bullet}="v"
    !{(1,1) }*{\bullet}="w"
    !{(-1,-1) }*{\bullet}="x"
    !{(0,0) }*{\otimes}="uv"
    !{(0,1) }*{\otimes}="uw"
    !{(0,-1) }*{\otimes}="vx"
    !{(1,0) }*{\otimes}="wv"
    !{(-1,0) }*{\otimes}="xu"
    "u"-"uw"-"w"-"wv"-"v"-"vx"-"x"-"xu"-"u"
    "u"-"uv"-"v"
  }
\]
  \caption{A graph $G$ and its edge subdivision $\sd G$.}
\end{figure}

In \cite{Lorenzini} Lorenzini first observed that the critical groups
$K(\sd G)$ and $K(G)$ determine each other in a trivial way, using the
description $K=Z^\#/Z$ as we now explain. If $\{C_1,\ldots,C_\beta\}$
is any set of directed cycles in $G$ whose incidence vectors
$\{z(C_i)\}_{i=1}^\beta$ give a $\ZZ$-basis for $Z_G$, then one can
subdivide those same cycles to obtain a $\ZZ$-basis $\{z_{\sd
  C_i}\}_{i=1}^\beta$ for $Z_{\sd G}$.  One then checks that
\[
\langle z_{\sd C_i}, z_{\sd C_j} \rangle = 2 \langle z(C_i), z(C_j) \rangle,
\]
for each $i,j$, since the inner product counts (with signs) the overlap
of edges between cycles $C_i, C_j$, and these overlaps double
in size after the subdivision.
Hence one has the following relation between their Gram
matrices:
\begin{equation}
\label{subdivision-Gram-matrix-relation}
M_{\sd G}^t M_{\sd G} = 2 M_G^t M_G
\end{equation}
and the following simple relation between their cokernels, the
critical groups:

\begin{proposition}[Lorenzini \cite{Lorenzini}]
\label{Lorenzini's-theorem}
Let $G$ be a graph with $\beta$ independent cycles.  Expressing $ K(G)
\cong \bigoplus_{i=1}^\beta \ZZ_{d_i} $ for positive integers
$d_1,d_2, \ldots,d_\beta \geq 1$, one has $ K(\sd G) \cong
\bigoplus_{i=1}^\beta \ZZ_{2d_i}$.
\end{proposition}
It will be useful later to have an expression of this result in terms of
explicit morphisms (as was done also in \cite{Lorenzini}).  Consider
the pair of adjoint maps defined $\RR$-linearly by
\[
  \begin{matrix}
    \RR^{E_{\sd G}} & \overset{h}{\longrightarrow} &\RR^{E_G} \\
    (u,uv)          & \longmapsto                  & (u,v) \\
    (uv,v)          & \longmapsto                  & (u,v) \\
                    &                              &       \\
    \RR^{E_G} & \stackrel{h^t}{\longrightarrow} &\RR^{E_{\sd G}} \\
    (u,v)          & \longmapsto                  & (u,uv)+(uv,v). 
  \end{matrix}
\]

One can easily check that these are morphisms of rational orthogonal
decompositions, and hence give rise to a morphism $h: K(\sd G)
\rightarrow K(G)$ of critical groups.  The relation
\eqref{subdivision-Gram-matrix-relation} between the two $\beta \times
\beta$ Gram matrices shows that the kernel-cokernel exact sequence
associated to $h$ takes this form:

$$
\begin{matrix}
0 \longrightarrow & \ker(h) 
  & \longrightarrow & K(\sd G) 
  & \overset{h}{\longrightarrow} & K(G)
  & \longrightarrow 0 \\
0 \longrightarrow & \ZZ_2^\beta
  & \longrightarrow & \bigoplus_{i=1}^\beta \ZZ_{2d_i} 
  & \longrightarrow & \bigoplus_{i=1}^\beta \ZZ_{d_i}
  & \longrightarrow 0. 
\end{matrix}
$$

Proposition~\ref{Lorenzini's-theorem} is equivalent to the assertion
that $K(\sd G)$ can be generated by $\beta$ elements and fits into an
exact sequence of this form, generalizing equation
\eqref{subdivision-complexity-relation} from the Introduction.

\subsection{A non-standard treatment of the complete graph}
\label{nonstandard-complete-graph-subsection}

Let $K_n$ be the complete graph on $n$ vertices.  A celebrated formula
of Cayley asserts that $\kappa(K_n)=n^{n-2}$ (see, e.g., \cite[Section
13.2]{Godsil}).  Generalizing this to compute the critical group
$K(K_n)$ is a favorite example of many papers in the subject.  We
approach this calculation in a slightly non-standard way here, mainly
because it will provide us with a crucial technical lemma for later
use in Section~\ref{kernel-subsection}.

\begin{proposition}
\label{complete-graph}
The complete graph $K_n$ has critical group
$$
K(K_n) \cong \ZZ_n^{n-2}.
$$
Furthermore, in the presentation $K(K_n) = \ZZ^E/(B \oplus Z)$, a
minimal generating set is provided by the images of any set of $n-2$
edges which form a spanning tree connecting $n-1$ out of the $n$
vertices.
\end{proposition}
\begin{proof}
  Since Cayley's formula implies $|K(K_n)| = |\ZZ_n^{n-2}|$, it will
  suffice to show that $K(K_n)$ is all $n$-torsion and that it can be
  generated by $n-2$ elements as in the second assertion.  Let
  $[n]:=\{1,2,\ldots,n\}$ denote the vertex set $V$ for $K_n$.

  To show $K(K_n)$ is all $n$-torsion, given any directed edge
  $e=(i,j)$ in $K_n$, we will prove that $n \cdot e$ is equal to a sum of
  cycles and bonds. Indeed, we can take the sum of the directed cycles $(i,j) + (j,k)
  + (k,i)$ for $k \in [n] - \{i,j\}$, and add the two bonds
  $$
  \begin{aligned}
    b(\{i\},[n]-\{i\}) &=(i,1) + (i,2) + \cdots + (i,n)\\
    b([n]-\{j\},\{j\}) &=(1,j) + (2,j) + \cdots + (n,j).
  \end{aligned}
  $$
  
  For the second assertion, let $T$ be a collection of $n-2$ edges
  that form a spanning tree connecting $n-1$ out of the $n$
  vertices. By symmetry, we may assume that $n$ is the vertex that is
  isolated by $T$. The edges of $K_n$ can be partitioned into two
  sets, $E(K_{n-1})$ and $\{(i,n)\}_{i=1}^{n-1}$.

  Any edge $e$ in $E(K_{n-1})$ either lies in $T$ or $T \cup \{e\}$
  contains a unique cycle
  that lets one express $e$ in terms of the elements of $T$ modulo
  $Z(K_{n-1})$, and hence modulo $Z=Z(K_{n})$.

  For each $1 \leq i \leq n-1$, the bond
  \[
  b_i := \sum_{\substack{ j=1 \\ j \neq i}}^n (i,j) \equiv 0 \mod B,
  \] 
  and it follows that
  \[ 
  (i,n) \equiv - \sum_{\substack{ j=1 \\ j \neq i}}^{n-1} (i,j) \mod B.
  \] 
  The edges in the sum all belong to $K_{n-1}$ and thus, according to
  the previous paragraph, can be written in terms of $T$ modulo
  $Z$. It follows that $(i,n)$ can be written in terms of $T$ modulo
  $B+Z$.\qedhere


\end{proof}

\section{Proof of Theorem~\ref{2-edge-connected-theorem}} 
\label{2-edge-connected-section}

Recall the statement of Theorem~\ref{2-edge-connected-theorem}:

\vskip .1in
\noindent
{\bf Theorem~\ref{2-edge-connected-theorem}.}
{\it
When the simple graph $G$ is $2$-edge-connected, the critical group $K(\line G)$ can
be generated by $\beta(G)$ elements.
}
\vskip .1in

The $\beta(G)$ generators will come from the set of edges in the complement $E \setminus T$ 
of a carefully chosen spanning tree $T$ for $G$.  For this we introduce
the following technical condition, which we have not encountered elsewhere.

\begin{definition} 
\label{absorption-order-definition}
\rm
For a connected graph $G=(V,E)$, say that a spanning tree $T \subset E$ for $G$ has an \emph{absorption order in} $G$ if one can
linearly order the union $V \sqcup T$ of its vertices and edges in the following way:
\begin{enumerate}
\item[(i)] The order begins with a vertex $v_0$ in $V$ followed by an edge $e_0$ of $T$, such that
$e_0$ is the unique edge of $T$ incident to $v_0$ (so $v_0$ is a leaf-vertex of $T$ attached along the leaf-edge $e_0$).
\item[(ii)]  For every other vertex $v$ in $V \setminus \{v_0\}$, there exists an edge $e=\{v,w\}$
such that $w$ occurs earlier in the order than $v$, and the edge $e$ either lies in $E \setminus T$ or
occurs earlier in the order than $v$.
\item[(iii)]  For every other edge $e$ in $T \setminus \{e_0\}$, there exists a vertex $v$
incident to $e$ which occurs earlier in the order than $e$,
and every other edge incident to $v$ either lies in $E \setminus T$
or occurs earlier in the order than $e$.
\end{enumerate}
\end{definition}

The relevance of an absorption order for a spanning tree is given by the algebraic consequence
in the following proposition.
Say that an orientation of the edges of a tree $T$ is {\it bipartite}
if, for every vertex $v$, the edges of $T$ incident to $v$ are either all oriented toward
$v$ or all oriented away from $v$.

\begin{proposition}
\label{absorption-order-significance}
Let $G=(V,E)$ be a simple graph, and assume it has a spanning tree $T \subset E$ which has an absorption order in $G$.

Then the images of the basis elements in $\ZZ^E$ corresponding to the edges $E \setminus T$ not lying on $T$
give a set of $\beta(G)$ generators for $K(\line G)$, using the presentation from Proposition~\ref{alternate-presentation}
$$
K(\line G) \cong \ZZ^E / \left(\ZZ e_0 + \partial_{\line G}(B_{\line G}) \right),
$$
assuming that the orientation chosen for $G$ restricts to a bipartite
orientation of $T$ (although $\line G$ may be oriented arbitrarily),
and the edge $e_0$ is the designated leaf-edge of $T$ appearing second
in the absorption order.

\end{proposition}

To prove this, note the following crucial lemma: 

\begin{lemma}
\label{bond-difference-lemma}
When a connected graph $G=(V,E)$ is oriented in a way that restricts to a bipartite orientation for
a spanning tree $T \subset E$, then any edge $e=\{v,w\}$ has 
$$
b_G(v) \equiv \pm b_G(w) \quad \mod \quad \ZZ e + \ZZ (E\setminus T) + \partial_{\line G} ( B_{\line G} ).
$$
\end{lemma}
\begin{proof}[Proof of Lemma~\ref{bond-difference-lemma}.]
Label the edges of $G$ incident to $v$ other than $e$ by
$$
\underbrace{e_1,\ldots,e_p,}_{\text{ in }T}
\underbrace{e_{p+1},e_{p+2},\ldots,e_P,}_{\text{ in }E \setminus T}
$$
and those incident to $w$ other than $e$ by
$$
\underbrace{f_1,\ldots,f_q,}_{\text{ in }T}
\underbrace{f_{q+1},f_{q+2},\ldots,f_Q}_{\text{ in }E \setminus T}.
$$
With these notations, one then has
\begin{equation}
\label{p-and-q-compare-1}
\begin{aligned}
\partial_{\line G} b_{\line G}(e) 
   &= (e_1-e)+\cdots+(e_P-e)+(f_1-e)+\cdots+(f_Q-e) \\
   &= (e_1+\cdots+e_P)+(f_1+\cdots+f_Q)-(P+Q)e.
\end{aligned}
\end{equation}
Because the orientation of $G$ when restricted to $T$ is bipartite, 
\begin{equation}
\label{p-and-q-compare-2}
\begin{aligned}
b_G(v) &= \pm( e_1+\cdots+e_p) \pm e_{p+1} \pm e_{p+2} \pm \cdots \pm e_P \\
b_G(w) &= \pm( f_1+\cdots+f_q) \pm f_{q+1} \pm f_{q+2} \pm \cdots \pm f_Q.
\end{aligned}
\end{equation}
Comparison of \eqref{p-and-q-compare-1} and \eqref{p-and-q-compare-2} shows that one of the two expressions
$b_G(v) + b_G(w)$ or $b_G(v) - b_G(w)$ differs from $\partial_{\line G} b_{\line G}(e)$ by
a $\ZZ$-linear combination of the edges in
$$
\{e\} \cup \{e_{p+1},e_{p+2},\cdots,e_P,f_{q+1},f_{q+2},\cdots,f_Q\}.
$$
Since the second set in the above union lies in $E \setminus T$, the lemma follows.
\end{proof}

\begin{proof}[Proof  of Proposition~\ref{absorption-order-significance}.]
  One needs to show that the subgroup of $\ZZ^E$ defined by
  $$
  I:=\ZZ (E \setminus T) + \ZZ e_0 + \partial_{\line G}( B_{\line G} )
  $$
  is {\it all} of $\ZZ^E$. Since $E \setminus T$ is a subset of $I$,
  it is enough to show that every edge $e$ in $T$ lies in $I$.  More
  strongly, one shows by induction on their location in the absorption
  order for $T$, that not only does every edge $e$ in $T$ lie in $I$,
  but also every vertex $v$ in $V$ has $b_G(v)$ lying in $I$.

  The base case for this induction deals both with the first vertex
  $v_0$ and the first edge $e_0$, which come at the beginning of the
  absorption order.  Since $v_0$ is a leaf vertex of $T$ along the
  edge $e_0$, one has $b_G(v_0)$ in $I$, since the only edges incident to
  $v_0$ are $e_0$ and edges of $E \setminus T$.  For the edge $e_0$,
  note that it trivially lies in $I$.

  In the inductive step, the next element in the absorption order is
  either a vertex $v \neq v_0$ or an edge $e \neq e_0$.

  If the next element is a vertex $v \neq v_0$, then by
  Definition~\ref{absorption-order-definition}(ii), there exists an
  edge $e=\{v,w\}$ for which $b_G(w)$ lies in $I$ by induction, and
  either $e$ lies in $E \setminus T$ (so that $e$ is in $I$) or $e$ is
  earlier in the order than $v$ (so that $e$ is in $I$ by induction).
  Hence Lemma~\ref{bond-difference-lemma} shows that $b_G(v)$ also lies in
  $I$.

  If the next element is an edge $e \neq e_0$, then by
  Definition~\ref{absorption-order-definition}(iii), there exists a
  vertex $v$ incident to $e$ for which $b_G(v)$ lies in $I$ by
  induction.  Note that $b_G(v)$ is a $\pm 1$ combination of all the
  edges $e'$ incident to $v$, and all of these other edges $e' \neq e$
  either have $e'$ in $E \setminus T$ (so that $e'$ is in $I$) or $e'$
  is earlier in the order than $e$ (so that $e'$ is in $I$ by
  induction).  Hence $e$ also lies in $I$.
\end{proof}

To show that $2$-edge-connected graphs $G$ always have a spanning tree
$T$ with an absorption order, we recall the well-known reformulation
of 2-edge-connectivity in terms of ear decompositions; see
e.g., \cite[Definition 4.2.7]{West}.

\begin{definition}\label{ear-def}
  Let $G=(V,E)$ be a simple graph.  An \it ear \rm of $G$ is a walk
  that alternates (incident) vertices $u_i$ and edges $e_i$
  \begin{equation}
    \label{ear-walk-labelling}
    v:=u_1,e_1,u_2,e_2,\ldots,u_\ell,e_\ell,u_{\ell+1}:=w
  \end{equation}
  such that the internal vertices $u_2,\ldots,u_{\ell}$ are each of
  degree $2$ in $G$.  If $v \neq w$, it is called an {\it open ear}
  (and necessarily $\ell \geq 1$), while if $v = w$, it is called a
  {\it closed ear} (and necessarily $\ell \geq 3$, because $G$ is
  simple).
  
  An \it{ear decomposition} \rm of $G$ is a decomposition of its
  vertices and edges
  $$
  P_0 \cup P_1 \cup \cdots \cup P_k
  $$ 
  such that $P_0$ is a cycle, and for $1 \leq i \leq k$, $P_i$ is an
  ear of $P_0 \cup P_1 \cup \cdots \cup P_i$.
\end{definition}
\begin{figure}
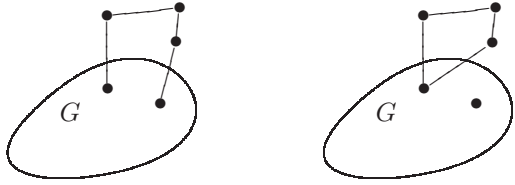

    \[ 
      \xygraph{
      !{<0cm,0cm>;<1cm,0cm>:<0cm,1cm>::}
      !~-{@{-}@[|(8.5)]}
      !{(0,0) }*{}="a"
      !{(1,1) }*{}="b"
      !{(2.5,0.5) }*{}="c"
      !{(2,-1)}*{}="d"
      !{(-1,-1)}*{}="e"
      !{(.5,0)}*{G}
      !{(1,.3)}*{\bullet}="u"
      !{(1.7,.1)}*{\bullet}="v"
      !{(-1,0);a(65)**{}?(1.4)}*{\bullet}="w"
      !{(-1,0);a(50)**{}?(1.8)}*{\bullet}="y"
      !{(-1,0);a(35)**{}?(1.6)}*{\bullet}="x"
      "a" -@`{"b", "c", "d","e"} "a"
      "u"-"w"-"y"-"x"-"v"
    } \qquad
    \xygraph{
      !{<0cm,0cm>;<1cm,0cm>:<0cm,1cm>::}
      !~-{@{-}@[|(8.5)]}
      !{(0,0) }*{}="a"
      !{(1,1) }*{}="b"
      !{(2.5,0.5) }*{}="c"
      !{(2,-1)}*{}="d"
      !{(-1,-1)}*{}="e"
      !{(.5,0)}*{G}
      !{(1,.3)}*{{}\bullet}="u"
      !{(1.7,.1)}*{\bullet}="v"
      !{(-1,0);a(65)**{}?(1.4)}*{\bullet}="w"
      !{(-1,0);a(50)**{}?(1.8)}*{\bullet}="y"
      !{(-1,0);a(35)**{}?(1.6)}*{\bullet}="x"
      "a" -@`{"b", "c", "d","e"} "a"
      "u"-"w"-"y"-"x"-"u"
    } 
    \]
    \caption{Graphs with open and closed ears.}
    \label{fig:eardecomp}
  \end{figure}
\begin{proposition}[{\cite[Theorem 4.2.8]{West}}] 
  \label{ear-prop}
  A graph is 2-edge-connected if and only if it has an ear
  decomposition.
\end{proposition}

In light of Proposition~\ref{absorption-order-significance}, the
following result implies Theorem~\ref{2-edge-connected-theorem}.
\begin{proposition}
\label{trees-with-absorption-orders-exist}
Let $G=(V,E)$ be a simple $2$-edge-connected graph.  Then $G$ has at
least one spanning tree $T \subset E$ with an absorption order in $G$.
\end{proposition}
\begin{proof}
  Induct on the number $k$ of ears in an ear decomposition $P_0 \cup
  P_1 \cup \cdots \cup P_k$ for $G$.

  In the base case $k=0$, the graph $G=P_0$ is a $n$-cycle.  Label its
  vertices $V=\{v_0,v_1,\ldots,v_{n-1}\}$ and edges
  $E=\{e_0,e_1,\ldots,e_{n-1}\}$ so that $e_i=\{v_i,v_{i+1}\}$ with
  indices taken modulo $n$.  Then one can easily check that
  $T=\{e_0,e_1,\ldots,e_{n-2}\}$ is a spanning tree, and
  $(v_0,e_0,v_1,e_1,\ldots,v_{n-2},e_{n-2},v_{n-1})$ is an absorption
  order for $T$ in $G$.

  In the inductive step, one may assume that $G^-:=P_0 \cup P_1 \cup
  \cdots \cup P_{k-1}$ has a spanning tree $T^-$ with an absorption
  order in $G^-$.  Choose the labelling of the endpoints $v,w$ of the
  ear $P_k$ so that $v$ comes weakly earlier than $w$ in the
  absorption order for $T^-$, where the vertices and edges of $P_k$
  are labelled as in \eqref{ear-walk-labelling}.  Extend $T^-$ to
  $$
  T:=T^- \sqcup \{u_2,u_3,\ldots,u_{\ell}\}
  $$ 
  which is easily seen to be a spanning tree for $G$.  One extends the
  absorption order for $T^-$ in $G^-$ to one for $T$ in $G$ by
  inserting the subsequence
  \begin{equation}
    \label{ear-subsequence}
    (u_2,e_2,u_3,e_3,u_4,\ldots,u_{\ell},e_{\ell})
  \end{equation}
  into the absorption sequence for $T^-$ in one of two possible locations,
  depending upon whether $v$ and $w$ are the initial vertex $v_0$ of
  the absorption order of $T^-$, or not.

  First we assume that $P_k$ is an open ear (that is, $v \neq w$) or
  $P_k$ is a closed ear with $v=w \neq v_0$. In this case, one can
  check that inserting the subsequence \eqref{ear-subsequence} {\it
    immediately after $v$} in the absorption order for $T^-$ in $G^-$
  gives an absorption order for $T$ in $G$.

  In the case that $P_k$ is a closed ear with $v=w=v_0$, one checks
  that inserting the subsequence \eqref{ear-subsequence} {\it at the
    very beginning} of the absorption order for $T^-$ in $G^-$ gives
  an absorption order for $T$ in $G$. Note that $u_2,e_2$ become the
  ``new'' $v_0, e_0$ in this absorption order. \qedhere
\end{proof}

We remark that the converse of 
Proposition~\ref{trees-with-absorption-orders-exist} is
false.  For example, one can check that 
the simple graph $G$ on vertex set $V=\{1,2,3,4\}$ with
edges $\{12,13,23,34\}$, which is {\it not} $2$-edge-connected,
{\it does} have an absorption ordering for any
choice of a spanning tree $T$ in $G$.
We have not investigated extensively the problem of
characterizing which graphs $G$
contain a spanning tree $T$ with an absorption ordering.

\section{Proof of Theorem~\ref{divisibility-theorem}}
\label{divisibility-section}
For a prime $p$ let $k(p)$ be the largest integer such that $p^{k(p)}$ divides
all of the sums $\deg_G(v)+\deg_G(w)$ as one runs through all edges
$e=\{v,w\}$ in the edge set $E$ of $G$. The goal of this section is to
give a proof of Theorem~\ref{divisibility-theorem}, which we now
recall.

\begin{divtheorem*}
  Let $G=(V,E)$ be a connected simple graph that contains at least one
  cycle of even length. Use the abbreviated notation $K:=K(\line
  G)$, and let $p$ be a prime for which the quantity $k(p) \geq 1$.
  
  Then for $G$ bipartite, one has
  $$
  K /p^{k(p)} K 
  \cong \ZZ_{p^{k(p)}}^{\beta(G)-1} \oplus \ZZ_{\gcd(p^k, |V|)},
  $$
  while for $G$ nonbipartite, one has
  $$
  K /p^{k(p)} K
  \cong \ZZ_{p^{k(p)}}^{\beta(G)-2} \oplus 
  \begin{cases}
    0 &\text{ if }p\text{ is odd,} \\
    \ZZ_2^2 &\text{ if }p=2\text{ and }|V|\text{ is even,} \\
    \ZZ_4 &\text{ if }p=2\text{ and }|V|\text{ is odd.}  
  \end{cases}
  $$  
\end{divtheorem*}
\begin{proof}
  One works again with the presentation from
  Proposition~\ref{alternate-presentation}
$$
K:=K(\line G) =\ZZ^E/\left(\ZZ e_0 + \ZZ\left( \partial_{\line G}
    b_{\line G}(e) \right)_{e \in E} \right)
$$
for some choice of an edge $e_0$ in $E$.  Given a vertex $v$ in $V$,
let $\delta_G(v)$ denote the element of $\ZZ^E$ which is the sum with
coefficient $+1$ of the basis elements in $\ZZ^E$ corresponding to
edges incident with $v$.  Given any edge $e=\{v,w\}$ in $E$, reasoning
as in equation \eqref{p-and-q-compare-1}, one finds that
$$
\partial_{\line G} b_{\line G}(e) 
 = \delta_G(v)+\delta_G(w)-(\deg_G(v)+\deg_G(w)) e.
$$
Letting $q:=p^{k(p)}$, one has therefore in $K/qK$ the relation
$$
\partial_{\line G} b_{\line G}(e) \equiv \delta_G(v)+\delta_G(w)
$$
and one can write a presentation for $K/qK$ as
\begin{equation}
\label{divisible-presentation}
K/qK=\ZZ_q^E / \left(\ZZ_q e_0 
  + \ZZ_q ( \delta_G(v)+\delta_G(w))_{e=\{v,w\} \in E} \right).
\end{equation}

We now make a particular choice of the edge $e_0$ for
this presentation, and exhibit a subset of $E$ having size
$\beta(G)-2$ or $\beta(G)-1$ which will
represent $\ZZ_q$-linearly independent elements in $K/qK$.
Because $G$ contains an even-length (not necessarily minimal) cycle,
it is possible to choose an edge $e_0$ in $E$ which lies on a {\it minimal}
cycle, so that $E \setminus \{e_0\}$ still connects all of $V$, and
so that $E \setminus \{e_0\}$ contains at least one odd cycle
in the case where $G$ is nonbipartite.  Now, in the bipartite case, pick 
$S \subset E \setminus \{e_0\}$ to be minimal with respect to the
property that $S$ connects all of $V$. In the non-bipartite case, pick $S$ to be minimal with respect to the following three properties: first, $S$ must connect all of $V$; second, $S$ must contain a unique cycle; and third, this cycle must be of odd length. This means that when $G$ is bipartite, $S$ is a spanning tree that avoids $e_0$, 
and when $G$ is nonbipartite, $S$ is a unicyclic spanning subgraph that 
avoids $e_0$, whose unique cycle $C$ is of odd length.

We first wish to show that, in either case, the images of the
elements $E \setminus S \setminus \{e_0\}$ are  $\ZZ_q$-linearly independent
in the presentation \eqref{divisible-presentation};  note that this
set $E\setminus S \setminus \{e_0\}$ has cardinality $\beta(G)-1$ when
$G$ is bipartite, and cardinality $\beta(G)-2$ when $G$ is nonbipartite.

So assume that $E \setminus S \setminus \{e_0\}$ are  $\ZZ_q$-linearly dependent
in $K/qK$.  Grouping the $\ZZ_q$-coefficients $c_v$ in front of each $\delta_G(v)$, one would have a sum
$\sum_{v \in V} c_v \delta_G(v)$ lying in 
$\ZZ_q e_0 + \ZZ_q(E\setminus S \setminus \{e_0\})$.
Thus this sum should have zero coefficient on every edge $e=\{v,w\}$ in $S$,
implying that $c_v=-c_w$ for every such edge.  Because $S$ is a spanning set
of edges, this forces the existence of a constant $c$ in $\ZZ_q$ for which
every $v$ in $V$ has $c_v = \pm c$.  In fact, when $G$ is bipartite with
vertex bipartition $V=V_1 \sqcup V_2$, this forces $c_{v_1}=c=-c_{v_2}$
for all $v_1$ in $V_1$ and $v_2$ in $V_2$, while for $G$ nonbipartite,
the existence of the odd cycle $C$ inside $S$ forces $c_v=c=-c$ for
all $v$ in $V$.  In either case, this means that $c_v=-c_w$ for all
edges $e=\{v,w\}$ in $E$, and hence the sum $\sum_{v \in V} c_v \delta_G(v)$
is actually zero in $\ZZ_q^E$.  Thus the linear independence is trivial.

It only remains now to analyze the quotient 
\begin{align}
\label{leftover-quotient}
\frac{K/qK}{\ZZ_q(E\setminus S \setminus \{e_0\})}
 = 
\ZZ_q^E / \left( \ZZ_q(E \setminus S) 
    + \ZZ_q (\delta_G(v)+\delta_G(w))_{e=\{v,w\} \in E} \right).  
\end{align}
Note that when $m$ is odd, for any sequence of vertices
$v_0,v_1,\ldots,v_{m-1},v_m$ one has a telescoping alternating sum
$$
\sum_{i=0}^{m-1} (-1)^i \left(\delta_G(v_i)+\delta_G(v_{i+1})\right) = \delta_G(v_0) + \delta_G(v_m).
$$
Also note that $S$ will contain paths of edges of odd length between
\begin{itemize}
\item every pair $(v_1,v_2)$ in $V_1 \times V_2$ when $G$ is bipartite, and
\item every ordered pair $(v,w)$ in $V \times V$ when $G$ is nonbipartite.
\end{itemize}
Thus, in either case, one has 
$$
\ZZ_q (\delta_G(v)+\delta_G(w))_{e=\{v,w\} \in E}
\quad = \quad 
\ZZ_q (\delta_G(v)+\delta_G(w))_{e=\{v,w\} \in S}.
$$ 
Using this last equation, one can rewrite the quotient on the right of
\eqref{leftover-quotient} as
\begin{equation}
\label{leftover-quotient-rewritten}
 \ZZ_q^S / \ZZ_q (\delta_S(v)+\delta_S(w))_{e=\{v,w\} \in S} \\
\end{equation}
where here we regard $S$ itself as a graph, namely the edge-induced
subgraph of $G=(V,E)$ having the same vertex set $V$ and edge set $S
\subset E$.

Note that this last expression in \eqref{leftover-quotient-rewritten}
does not depend upon the ambient graph $G$, but only on the subgraph
$S$.  We therefore rename it $K_q(S)$ to emphasize this dependence on
$S$ alone.  It remains to analyze this group $K_q(S)$ in both the
bipartite and nonbipartite cases.

\noindent \emph{Case 1}: $G$ is bipartite (and hence so is $S$). It follows
that $S$ is a spanning tree on $V$, with vertex bipartition $V = V_1
\sqcup V_2$.  By the above discussion,
$$
\begin{aligned}
K_q(S) 
 &= \ZZ_q^S / \ZZ_q ( \delta_S(v)+\delta_S(w) )_{(v_1,v_2) \in V_1 \times V_2} \\
 &= \ZZ_q^S / \ZZ_q \left( \sum_{v \in V} c_v \delta_S(v): 
              \sum_{v_1 \in V_1} c_{v_1} = \sum_{v_2 \in V_2} c_{v_2} 
                \text{ in }\ZZ_q \right).
\end{aligned}
$$

We first show by induction on $|S|$ that $K_q(S)$ is cyclic, generated by the image of
any {\it leaf edge} $e$ of $S$, that is, an edge $e$ incident to some {\it leaf vertex} $v$
having $\deg_S(v)=1$.  The base case $|S|=1$ is trivial.  In the inductive
step, pick another leaf edge $e'$ in $S$; we will show it has image $0$ in 
the quotient $K_q(S)/\ZZ_q e$.  If $e'$ is incident to leaf vertex $v'$, then 
for any $c$ in $\ZZ_q$, one has 
$$
e' + c e = \delta_S(v') + c \delta_S(v).
$$
Taking $c=-1$ (respectively $+1$) when $v,v'$ lie in the same
(resp. different) set $V_1$ or $V_2$, one obtains an element that is
zero in $K_q(S)$, and hence $e' \equiv 0$ in $K_q(S)/\ZZ_q e$.  Now, replacing $S$ by $S \setminus \{e'\}$,
one can induct on $|S|$,
completing the inductive step and showing that $K_q(S)$ is generated
by $e$.

We next analyze the order of this cyclic generator $e$ within
$K_q(S)$.  We claim that $c \cdot e = 0$ in $K_q(S)$ if and only if
$c$ lies in $|V| \ZZ_q$.  This would finish the proof in the bipartite
case, as it would show that $K_q(S)$ is isomorphic to the subgroup of
$\ZZ_q$ generated by the element $|V|$.  This subgroup is isomorphic
to $\ZZ_{\gcd(q,|V|)}$, where $q=p^k$.  Hence this would imply $K/qK
\cong \ZZ_{p^{k(p)}}^{\beta(G)-1} \oplus \ZZ_{\gcd(p^{k(p)},|V|)}$, as
desired.

To see the claim, assume that $c \cdot e = 0$ in $K_q(S)$ for some $c$
in $\ZZ_q$.  This means one has a sum $\sum_{v \in V} c_v \delta_G(v)
=c \cdot e$ in which $\sum_{v_1 \in V_1} c_{v_1} = \sum_{v_2 \in V_2}
c_{v_2}$.  This happens if and only if the sum has zero coefficient on
all edges $e'$ in $S \setminus \{e\}$.  If $e=\{v,w\}$ with the leaf
vertex $v$ lying in $V_1$, and $w$ in $V_2$, this means
$c_{v_2}=c_w=-c_{v_1}$ for all $v_1 \in V_1 \setminus \{v\}$ and $v_2
\in V_2 \setminus \{w\}$.  Then the condition $\sum_{v_1 \in V_1}
c_{v_1} = \sum_{v_2 \in V_2} c_{v_2}$ forces
$$
c_v+(|V_1|-1)(-c_w)=|V_2|(c_w)
$$ 
i.e., $c_v=(|V|-1)c_w$.  Hence this can occur if and only if
$c=c_v+c_w=|V|c_w$, that is, if $c$ lies in $|V| \ZZ_q$.

\noindent \emph{Case 2}. $G$ is nonbipartite (and hence so is $S$). In this
case $S$ is a spanning unicyclic graph, whose unique cycle $C$ is of
odd length.  By the above discussion,
$$
\begin{aligned}
K_q(S) 
 &= \ZZ_q^S / \ZZ_q ( \delta_S(v)+\delta_S(w) )_{(v,w) \in V \times V} \\
 &= \ZZ_q^S / 
     \left\{ \sum_{v \in V} c_v \delta_S(v): 
     \sum_{v \in V} c_{v} \in 2\ZZ_q \right\}.
\end{aligned}
$$
Thus, if one defines the tower of $\ZZ$-lattices (i.e., free abelian
groups)
$$
L:=\ZZ^S \ \supset\ 
   M:=\ZZ(\delta_S(v))_{v \in V}  \ \supset \  
    N:=\left\{ \sum_{v \in V} c_v \delta_S(v): 
                      \sum_{v \in V} c_{v} \in 2\ZZ\right\},
$$
then one has a short exact sequence
\begin{equation}
\label{lifted-sequence}
0 \rightarrow  \underbrace{\frac{M\otimes_\ZZ \ZZ_q}{N\otimes_\ZZ \ZZ_q}}_{(M/N) \otimes_\ZZ \ZZ_q}
\rightarrow  \underbrace{\frac{L\otimes_\ZZ \ZZ_q}{N\otimes_\ZZ \ZZ_q}}_{K_q(S)}
\rightarrow  \underbrace{\frac{L\otimes_\ZZ \ZZ_q}{M\otimes_\ZZ \ZZ_q}}_{(L/M) \otimes_\ZZ \ZZ_q}  \rightarrow 0
\end{equation}
Here we have used on the two ends of the sequence the fact\footnote{That is, taking tensor products $(-)_{\otimes \ZZ} \ZZ_q$ is
  {\it right exact}, so when applied to the exact sequence
$$B \rightarrow A \rightarrow A/B \rightarrow 0$$
it gives the exact sequence
$$
B \otimes_\ZZ \ZZ_q \rightarrow A \otimes_\ZZ \ZZ_q \rightarrow (A/B)
\otimes_\ZZ \ZZ_q \rightarrow 0.
$$
} that for any pair of nested abelian groups $B \subset A$, one has
$$
\left( A \otimes_\ZZ \ZZ_q \right)/ \left( B \otimes_\ZZ \ZZ_q \right)
\cong
\left( A/B \right) \otimes_\ZZ \ZZ_q.
$$
Furthermore, it is easy to see that both $M/N$ and $L/M$ are
isomorphic to $\ZZ_2$:
\begin{itemize}
\item
The isomorphism $M/N \cong \ZZ_2$ comes from choosing any $\ZZ$-basis for the
lattice $M$, thus identifying $M \cong \ZZ^{|V|}$, and noting that
under this identification, $N$ is identified with the index $2$ sublattice
$\{ x \in \ZZ^{|V|}: \sum_{v \in V} x_v \in 2\ZZ \}$.
\item
The isomorphism $L/M \cong \ZZ_2$ is equivalent to the assertion
that the square {\it (unsigned) edge-node incidence} matrix having columns indexed by 
the nodes $V$ and rows indexed by the edges $S$ will have determinant
$\pm 2$.  This is a well-known fact for connected unicyclic graphs $S$ whose
unique cycle $C$ is odd; see, e.g., \cite[p. 560, proof of Thm. 3.3]{Stanley-zonotopes}.
It is easily proven by first checking that the determinant is
scaled by $\pm 1$ when one removes a row and column corresponding to
a leaf edge and its incident leaf vertex in $S$.
This reduces the assertion to the case where $S=C$ is just an odd cycle
itself, where the determinant can be calculated directly via Laplace expansion.
\end{itemize}

Hence both of the outer terms $(M/N) \otimes_\ZZ \ZZ_q , (L/M)
\otimes_\ZZ \ZZ_q $ in the short exact sequence
\eqref{lifted-sequence} are isomorphic to $\ZZ_2 \otimes_\ZZ \ZZ_q$,
which vanishes for $p$ odd and equals $\ZZ_2$ for $p=2$.  Thus
\eqref{lifted-sequence} shows that $K_q(S)$ vanishes for $p$ odd, and
shows for $p=2$ that $K_q(S)$ is either $\ZZ_2^2$ or $\ZZ_4$.  To
distinguish these possibilities when $p=2$, we analyze the additive
orders of each edge $e$ in $S$ as elements of $K_q(S)$.

Note that for any leaf edge $e$ in $S$, say with leaf vertex $v$, one
has $e=\delta_S(v)$, and hence $2e=\delta_S(v)\equiv 0$ in $K_q(S)$.
Thus using a leaf-induction, one sees that any edge $e$ in $S
\setminus C$ has $2e \equiv 0$ in $K_q(S)$.

Meanwhile, we claim that for any edge $e=\{v,w\}$ in $C$, one has $c
\cdot e \equiv 0$ in $K_q(S)$ if and only if $|V| \cdot c$ lies in $4
\ZZ_q$.  To see the claim, assume that $c \cdot e = 0$ in $K_q(S)$ for
some $c$ in $\ZZ_q$.  This means one has a sum $\sum_{v \in V} c_v
\delta_G(v) =c \cdot e$ with $\sum_{v \in V} c_v \in 2\ZZ_q$.  This
happens if and only if the sum has zero coefficient on all edges $e'$
in $S \setminus \{e\}$.  Applying this for the edges $e'$ in $C
\setminus \{e\}$, one concludes that $c_v=c_w$, and hence
$c=c_v+c_w=2c_v$.  Applying this for the remaining edges $e'$ in $S
\setminus C$, one concludes that $c_w = \pm c_v$ for all $w$ in $V$.
But then the condition that $\sum_{w \in V} c_w$ lies in $2\ZZ_q$
means that $|V|\cdot c_v$ also lies in $2\ZZ_q$, i.e., that $|V|\cdot
c=2|V|\cdot c_v$ lies in $4\ZZ_q$.  One concludes that edges $e$ in
$C$ have order $2$ when $|V|$ is even, and order $4$ when $|V|$ is
odd.  Since every edge $e$ in $S \setminus C$ has $2e \equiv 0$ in $K_q(S)$, this implies $K_q(S)
\cong \ZZ_2^2$ when $|V|$ is even and $K_q(S) \cong \ZZ_4$ when $|V|$
is odd.
\end{proof}

\section{Proof of Theorem~\ref{regular-result}}
\label{regular-result-section}

Recall here the statement of Theorem~\ref{regular-result}.
\begin{regtheorem*}
  For any connected $d$-regular simple graph $G$ with $d \geq 3$,
  there is a group homomorphism 
  $$
  K(\line G) \overset{f}{\rightarrow} K(\sd G)
  $$
  whose kernel-cokernel exact sequence takes the form
  $$
  0 
  \rightarrow \ZZ_d^{\beta(G)-2} \oplus C
  \rightarrow K(\line G)
  \overset{f}{\rightarrow} K(\sd G)
  \rightarrow C
  \rightarrow 
  0
  $$
  in which the cokernel $C$ is the following cyclic $d$-torsion group:
  $$
  C=
  \begin{cases}
    0        & \text{ if } G \text{ is non-bipartite and }d\text{ is odd},\\
    \ZZ_2 & \text{ if } G \text{ is non-bipartite and }d \text{ is even},\\
    \ZZ_d & \text{ if } G \text{ is bipartite}.
  \end{cases}
  $$
\end{regtheorem*}

\subsection{Defining the morphism $f$}
\label{defining-f-subsection}

We begin our proof of the theorem by first defining a linear map
$f:\RR^{E_{\line G}} \rightarrow \RR^{E_{\sd G}}$ which will turn out
to be a morphism of rational orthogonal decompostions.

\begin{figure}
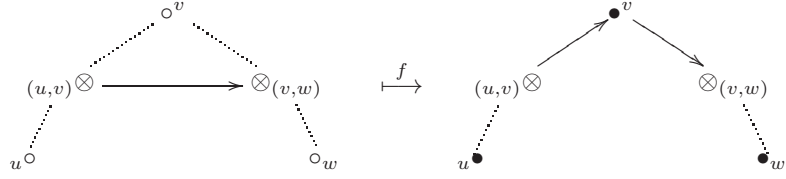

\[\xygraph{
!{<0cm,0cm>;<1cm,0cm>:<0cm,1cm>::}
!{(0,-1) }*+{{}_{u}\circ}="u" 
!{(.5,0) }*+{{}_{(u,v)}\otimes}="uv"
!{(2,1) }*+{\circ^{v}}="v"
!{(3.5,0)}*+{\otimes_{(v,w)}}="vw"
!{(4,-1)}*+{\circ_{w}}="w"
"u"-@{..}"uv"-@{..}"v"-@{..}"vw"-@{..}"w" 
"uv":"vw" 
}
\quad \stackrel{f}{\longmapsto}\quad
\xygraph{
!{<0cm,0cm>;<1cm,0cm>:<0cm,1cm>::}
!{(0,-1) }*{{}_{u}\bullet}="u" 
!{(.5,0) }*+{{}_{(u,v)}\otimes}="uv"
!{(2,1) }*{\bullet^{v}}="v"
!{(3.5,0)}*+{\otimes_{(v,w)}}="vw"
!{(4,-1)}*{\bullet_{w}}="w"
"u"-@{..}"uv":"v":"vw"-@{..}"w" 
}
\]
\caption{The action of $f$ on a single edge of $\line G$.}\label{fig:f-def}
\end{figure}

\begin{definition}\label{f-definition}
  Define a $\RR$-linear map $f:\RR^{E_{\line G}} \rightarrow
  \RR^{E_{\sd G}}$ by setting
  $$
  f(uv,vw) = (uv,v) + (v,vw)
  $$
  for every pair of edges $\{u,v\}, \{v,w\}$ of $G$ incident at some
  vertex $v$ (see Figure~\ref{fig:f-def}).  Equivalently, the adjoint
  map $f^t$ is defined by
  $$
  f^t(uv,v) = \sum_{w \in V: \{v,w\} \in E} (uv,vw).
  $$
\end{definition}

The following definitions and lemma will be useful both for showing that
$f$ gives a morphism, and in our later analysis.

\begin{definition}\label{local-global-defn}
Given a directed cycle 
$$
C=\{(v_1,v_2), (v_2,v_3), \ldots, (v_{m-1},v_m), (v_m,v_1) \}
$$
in $G$, let 
$$
\begin{aligned}
\sd C &:= \{(v_1,{v_1 v_2}),({v_1 v_2},v_2), 
            (v_2,{v_2 v_3}),({v_2 v_3},v_3), \ldots \}\\
\line C &:= \{(v_1 v_2, v_2 v_3) ,(v_2 v_3, v_3 v_4),  \ldots, 
              (v_{m-1}v_m, v_m v_1), (v_m v_1, v_1 v_2) \}
\end{aligned}
$$
denote corresponding cycles in $\sd G, \line G$.

Cycles in $\line G$ of the form $\line C$ where $C$ is a cycle of $G$
will be called {\it global} cycles.  A cycle in $\line G$ will be
called {\it local (to vertex $v$)} if every vertex $v_iv_j$ of $\line
G$ visited by the cycle has $v \in \{v_i, v_j\}$.

\end{definition}

\begin{lemma}
\label{induced-cycles}
Let $G$ be a graph, and let $\{C\}$ be a set of directed cycles
indexing a spanning set $\{z(C)\}$ for the cycle space $Z_G$. Then
\begin{enumerate}
\item $Z_{\sd G}$ will be spanned by the incidence vectors
  $\{z(\sd C) \}$ of the associated subdivided cycles, and
\item $Z_{\line G}$ will be spanned by the incidence vectors $\{ z(\line C) \}$ for
their associated global cycles together with all local cycles.
\end{enumerate}
\end{lemma}
\begin{proof}
  Assertion (1) of was implicitly used in
  Section~\ref{subdivision-example}, and should be clear either from
  elementary algebraic topology or from the discussion of bases for
  $Z_G$ coming from spanning forests at the beginning of
  Section~\ref{cycles-bonds-section}.

  For assertion (2), given any directed cycle in $\line G$, put an
  equivalence relation on its edges by taking the transitive closure
  of the following relation: two consecutive edges $(uv,vw), (vw,wx)$
  in the cycle are equivalent if there exists a vertex $y$ of $G$
  contained in $\{u,v\} \cap \{v,w\} \cap \{w,x\}$.  The global cycles
  in $\line G$ are by definition those in which the equivalence
  classes for this relation all have cardinality two (N.B.: here one
  is using the assumption that $G$ is simple).  Given a cycle $z$ in
  $\line G$ that contains equivalence classes of size at least $3$, it
  is easy to see that one can always add a local cycle to $z$ and
  reduce the number of such equivalence classes: if the equivalence
  class and its neighbors in $z$ correspond to these terms
  $$
  \cdots + (ab_1,yb_1)+(yb_1,yb_2)+(yb_2,yb_3)+\cdots
  +(yb_{t-1},yb_t)+(yb_t,b_tc)+\cdots
  $$
  where $a,c \neq y$, then subtracting the local cycle
  $$
  (yb_1,yb_2)+(yb_2,yb_3)+\cdots +(yb_{t-1},yb_t)+(yb_t,yb_1)
  $$
  gives a result that looks locally like
  \[
  \cdots + (ab_1,yb_1) + (yb_1,yb_t) + (yb_t,b_tc)+ \cdots . \qedhere
  \]
\end{proof}

\begin{corollary}
  For any $d$-regular simple graph $G$, the map $f: \RR^{E_{\line G}}
  \rightarrow \RR^{E_{\sd G}}$ from Definition~\ref{f-definition} is a
  morphism of the associated rational orthogonal decompositions, and
  hence induces a group homomorphism
  $$ 
  f: K(\line G) \rightarrow K(\sd G).
  $$
\end{corollary}

\begin{figure}
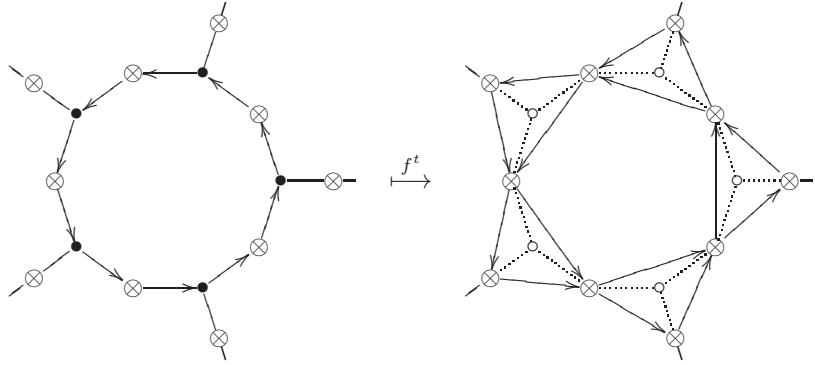

\[
\xygraph{
!{<0cm,0cm>;<1cm,0cm>:<0cm,1cm>::}
!{(0,0);a(0)**{}?(1.5)}*{\bullet}="a0"
!{(0,0);a(0)**{}?(2.2)}*{\otimes}="aa0"
!{(0,0);a(0)**{}?(2.5)}*{}="aaa0"
!{(0,0);a(36)**{}?(1.5)}*{\otimes}="a1"
!{(0,0);a(72)**{}?(1.5)}*{\bullet}="a2"
!{(0,0);a(72)**{}?(2.2)}*{\otimes}="aa2"
!{(0,0);a(72)**{}?(2.5)}*{}="aaa2"
!{(0,0);a(108)**{}?(1.5)}*{\otimes}="a3"
!{(0,0);a(144)**{}?(1.5)}*{\bullet}="a4"
!{(0,0);a(144)**{}?(2.2)}*{\otimes}="aa4"
!{(0,0);a(144)**{}?(2.6)}*{}="aaa4"
!{(0,0);a(180)**{}?(1.5)}*{\otimes}="a5"
!{(0,0);a(216)**{}?(1.5)}*{\bullet}="a6"
!{(0,0);a(216)**{}?(2.2)}*{\otimes}="aa6"
!{(0,0);a(216)**{}?(2.6)}*{}="aaa6"
!{(0,0);a(252)**{}?(1.5)}*{\otimes}="a7"
!{(0,0);a(288)**{}?(1.5)}*{\bullet}="a8"
!{(0,0);a(288)**{}?(2.2)}*{\otimes}="aa8"
!{(0,0);a(288)**{}?(2.5)}*{}="aaa8"
!{(0,0);a(324)**{}?(1.5)}*{\otimes}="a9"
"a0":"a1":"a2":"a3":"a4":"a5":"a6":"a7":"a8":"a9":"a0"
"a0"-"aa0"-"aaa0"
"a2"-"aa2"-"aaa2"
"a4"-"aa4"-"aaa4"
"a6"-"aa6"-"aaa6"
"a8"-"aa8"-"aaa8"
}
\quad\stackrel{f^t}{\longmapsto}\quad
\xygraph{
!{<0cm,0cm>;<1cm,0cm>:<0cm,1cm>::}
!{(0,0);a(0)**{}?(1.5)}*{\circ}="a0"
!{(0,0);a(0)**{}?(2.2)}*{\otimes}="aa0"
!{(0,0);a(0)**{}?(2.5)}*{}="aaa0"
!{(0,0);a(36)**{}?(1.5)}*{\otimes}="a1"
!{(0,0);a(72)**{}?(1.5)}*{\circ}="a2"
!{(0,0);a(72)**{}?(2.2)}*{\otimes}="aa2"
!{(0,0);a(72)**{}?(2.5)}*{}="aaa2"
!{(0,0);a(108)**{}?(1.5)}*{\otimes}="a3"
!{(0,0);a(144)**{}?(1.5)}*{\circ}="a4"
!{(0,0);a(144)**{}?(2.2)}*{\otimes}="aa4"
!{(0,0);a(144)**{}?(2.6)}*{}="aaa4"
!{(0,0);a(180)**{}?(1.5)}*{\otimes}="a5"
!{(0,0);a(216)**{}?(1.5)}*{\circ}="a6"
!{(0,0);a(216)**{}?(2.2)}*{\otimes}="aa6"
!{(0,0);a(216)**{}?(2.6)}*{}="aaa6"
!{(0,0);a(252)**{}?(1.5)}*{\otimes}="a7"
!{(0,0);a(288)**{}?(1.5)}*{\circ}="a8"
!{(0,0);a(288)**{}?(2.2)}*{\otimes}="aa8"
!{(0,0);a(288)**{}?(2.5)}*{}="aaa8"
!{(0,0);a(324)**{}?(1.5)}*{\otimes}="a9"
"a0"-@{..}"a1"-@{..}"a2"-@{..}"a3"-@{..}"a4"-@{..}"a5"-@{..}"a6"-@{..}"a7"-@{..}"a8"-@{..}"a9"-@{..}"a0"
"a0"-@{..}"aa0"-"aaa0"
"a2"-@{..}"aa2"-"aaa2"
"a4"-@{..}"aa4"-"aaa4"
"a6"-@{..}"aa6"-"aaa6"
"a8"-@{..}"aa8"-"aaa8"
"aa0":"a1":"aa2":"a3":"aa4":"a5":"aa6":"a7":"aa8":"a9":"aa0"
"a1":"a3":"a5":"a7":"a9":"a1"
}
\]
\caption{An example of a subdivided cycle in $\sd G$, and
  its image under $f^t$ in $\line G$ when $G$ is $3$-regular.}
\label{f-transpose}
\end{figure}

\begin{proof}
  By Proposition~\ref{Treumann-characterization}, one must show both
  $f(Z_{\line G}) \subset Z_{\sd G}$ and $f^t(Z_{\sd G}) \subset
  Z_{\line G}$.

  To show $f(Z_{\line G}) \subset Z_{\sd G}$, using
  Lemma~\ref{induced-cycles}(ii), it suffices to show that $f$ takes
  both global and local cycles in $\line G$ to cycles in $Z_{\sd G}$.
  This is easy (and requires no assumption about the $d$-regularity of
  $G$): local cycles map to $0$ under $f$, and a global cycle of the
  form $\line C$ satisfies $f(z(\line C)) = z(\sd C)$.

  To show $f^t(Z_{\sd G}) \subset Z_{\line G}$, using
  Lemma~\ref{induced-cycles}(i), it suffices to show for every
  directed cycle $C$ in $G$ that $f$ takes the subdivided cycle
  \[
  z(\sd C)=(v_1 v_2,v_2)+(v_2,v_2v_3)+(v_2v_3,v_3)+(v_3,v_3v_4)+
  \cdots+(v_kv_1,v_1)+(v_1,v_1v_2)
  \]
  to a sum of cycles in $Z_{\line G}$.  The regularity of $G$ implies
  that each $v_i$ has $d-2$ neighbors off the cycle; label them
  $u_i^1,\ldots,u_i^{d-2}$. Then one can write
  $$
  (f^t)(z(\sd C)) = 2 z(\line C) + \zeta_1 + \cdots + \zeta_{d-2}
  $$
  where for $j=1,2,\ldots,d-2$ one defines the element of $Z_{\line
    G}$
  \begin{align*}
  \zeta_j :=
  (v_1v_2,v_2u_2^j)+(v_2u_2^j,v_2v_3)+(v_2v_3,v_3u_3^j)+&(v_3u_3^j,v_3v_4)+
  \cdots \\+&(v_kv_1,v_1u_1^j)+(v_1u_1^j,v_1v_2).
  \end{align*}
  An example with $d=3$ is shown in Figure \ref{f-transpose},
  depicting the subdivided cycle $\sd C$ in $\sd G$, and then its
  image under $f^t$ in $\line G$, which decomposes into $2$ copies of
  the inner cycle $\line C$ along with $1$ $(=d-2)$ outer cycle
  $\zeta_1$.
\end{proof}

\subsection{The kernel and cokernel of $f$ are $d$-torsion}
\label{d-torsion-subsection}

\begin{proposition}
\label{f-scaling-prop}
For any $d$-regular connected graph $G$, both maps 
$$
\begin{aligned}
f^t f &: K(\line G) \rightarrow K(\line G)\\
f f^t &: K(\sd G) \rightarrow K(\sd G)
\end{aligned}
$$
are scalar multiplications by $d$.
\end{proposition}
\begin{proof}
  The proofs of these are straightforward computations:
\begin{align*}
  f^t f (uv,vw) &= f^t(uv,v) + f^t(v,vw) \\
  &= \sum_{x \in V:\{v,x\} \in E} (uv,vx) + (xv,vw) \\
  &= d \cdot (uv,vw) +
  \sum_{x \in V:\{v,x\} \in E} \left( (uv,vx) + (xv,vw) + (vw,uv) \right) \\
  &= d \cdot (uv,vw) \text{ mod } Z_{\line G} .\\
  f f^t (uv,v)  &= \sum_{x \in V:\{v,x\} \in E} f(uv,vx) \\
  &= \sum_{x \in V:\{v,x\} \in E} (uv,v) + (v,vx) \\
  &= d \cdot (uv,v) + \sum_{x \in V:\{v,x\} \in E} (v,vx) \\
  &=d \cdot (uv,v) \text{ mod } B_{\sd G}.\qedhere
\end{align*}
\end{proof}

\begin{corollary}
\label{f-torsion-corollary}
For any $d$-regular connected graph $G$, 
both $\ker(f)$ and $\coker(f)$ are all $d$-torsion.
\end{corollary}
\begin{proof}
For $x \in \ker(f)$ and $y \in \coker(f)$, one has 
\begin{align*}
d \cdot x &= f^t f(x) = f^t ( 0 ) = 0, \\
d \cdot y &= f f^t(y) \in \im(f). \qedhere
\end{align*}
\end{proof}

\subsection{Analyzing the cokernel}
\label{cokernel-subsection}
\begin{proposition}
\label{cokernel-description}
For any $d$-regular connected graph $G$, the group $C:=\coker(f)$ is a
cyclic group as described in Theorem~\ref{regular-result}.
\end{proposition}

\begin{proof}
  We will use the presentation
  \begin{equation}
    \label{coker-presentation}
    C:=\coker(f) : = K( \sd G ) / \im(f) 
    = \ZZ^{E_{\sd G}} / \left( B_{\sd G} + Z_{\sd G} + \im(f) \right),
  \end{equation}
  which follows from our first definition of the critical group
  (as in Section~\ref{rods-section}).

  To see that $C$ is cyclic, note that there are two ways for a pair
  of edges in $\sd G$ to be incident at a vertex, and in either case
  their images in $C$ will differ by a sign:
    \begin{align*}
    (u,uv) &= -(uv,v) \mod  B_{\sd G},\\
    (uv,v) &= -(v,vw) \mod \im(f).      
    \end{align*}

  Since $G$ is connected, this shows $C$ is cyclic, generated by the
  image of any directed edge of $\sd G$.  Furthermore, it is a
  quotient of $\ZZ_d$ by Corollary~\ref{f-torsion-corollary}.

  When $G$ is bipartite, in order to show $C = \ZZ_d$, it will suffice
  to exhibit a surjection $C \twoheadrightarrow \ZZ_d$.  Let the
  vertex set $V$ for $G$ have bipartition $V=V_1 \sqcup V_2$ and
  assume that all the edges of $G$ are oriented from $V_1$ to
  $V_2$. Define an abelian group homomorphism 
  \[
  \phi:\ZZ^{E_{\sd G}} \to \ZZ, \qquad \phi(v_1, v_1v_2) =
  \phi(v_1v_2,v_2) = 1,
  \]
  where $v_i \in V_i$ for $=1,2$. One can check that each of the three
  subgroups $B_{\sd G}, Z_{\sd G}, \im(f)$ by which one mods out in
  \eqref{coker-presentation} is mapped via $\phi$ into the subgroup
  $d\ZZ$:
  \begin{itemize}
  \item Any directed cycle $C$ in $\sd G$ has
    $\phi(z(C))=0$ (due to the fact that $C$ will have even length),
  \item any edge $e$ of $\line G$ has $\phi(f(e))=0$,
  \item any vertex $v_1v_2$ in $\sd G$ has $\phi(b_{\sd G}(v_1v_2))=0$,
  \item any vertex $v_i$ in $\sd G$ has $\phi(b_{\sd
      G}(v_i))= (-1)^{i-1} d$, where $i=1,2$.
  \end{itemize}
  Thus $\phi$ induces a surjection from $C$ onto $\ZZ_d$, as desired.

  If $G$ is not bipartite, it contains some (directed) odd cycle
  $C$. Pick any directed edge $e$ in the subdivision $\sd C$ and use
  the two relations (a), (b) to rewrite it successively as $\pm$ the
  other directed edges in the cycle.  It changes sign each time one
  uses (a) to pass through a vertex of $\sd C$ that comes from an edge
  of $G$.  Since there are an odd number of such edges in the cycle,
  it will change sign an odd number of times before it returns,
  yielding
  $$
  e = -e \text{ mod } B_{\sd G} + \im(f).
  $$
  Hence $2e =0$ in $C$, so $C$ is a quotient of $\ZZ_2$.
  
  Since $C$ is also a quotient of $\ZZ_d$, when $d$ is odd, one must
  have $C=0$.  When $d$ is even, consider the index $2$ sublattice
  $\Lambda$ of $\ZZ^{E_{\sd G}}$ consisting of those vectors whose sum
  of coordinates is even.  Without any parity assumption on $d$, it is
  true that $\im(f) \subset \Lambda$ (by definition of $f$) and
  $Z_{\sd G} \subset \Lambda$ (because the subdivided cycles $\sd C$
  have evenly many edges).  The assumption that $d$ is even implies
  that $B_{\sd G}$ also lies in $\Lambda$: $B_{\sd G}$ is generated by
  the bonds in $\sd G$ of the form $b_{\sd G}(v)$ for vertices $v$
  of $\sd G$, and every vertex in $\sd G$ has degree either $2$ or
  $d$.  Consequently, the presentation \eqref{coker-presentation}
  shows that $C$ surjects onto $\ZZ^{E_{\sd g}}/\Lambda \cong \ZZ_2$.
\end{proof}

\subsection{Analyzing the kernel}
\label{kernel-subsection}

It remains to understand $\ker(f)$, or equivalently by
Proposition~\ref{identify-kernel-cokernel}, to understand its
Pontrjagin dual
\begin{equation}
\label{transpose-cokernel-presentation}
\coker(f^t) = \ZZ^{E_{\line G}}/ \left( Z_{\line G} + B_{\line G} + \im(f^t) \right).
\end{equation}
This will come about by reformulating this presentation, in order to
analyze it {\it locally}.

\begin{definition}
  For each vertex $v \in V_G$ of a $d$-regular simple graph
  $G=(V_G,E_G)$, define inside $\line G$ the {\it $d$-clique local to $v$}
  $$
  K^{(v)}_d=(V(K^{(v)}_d), E(K^{(v)}_d))
  $$ 
  to be the vertex-induced subgraph of $\line G$ on the vertex set
  $$
  V(K^{(v)}_d) :=\{vw: vw \in E_G\}.
  $$  
\end{definition}

Note that the edges of $\line G$ form a disjoint decomposition
\begin{equation}
\label{disjoint-edge-decomposition}
E_{\line G}=\bigsqcup_{v \in V_G} E(K^{(v)}_d)
\end{equation}
since $G$ was assumed to be a simple graph.  Also note that a cycle in
$\line G$ is local to vertex $v$, as in
Definition~\ref{local-global-defn}, if and only if it is supported on
the edges $ E(K^{(v)}_d)$.  If one lets $Z^{global}_{\line G}$ be the
span of global cycles $\{z_{\line C}\}$ coming from any spanning set
of cycles $\{z_C\}$ for $Z_G$, then Lemma~\ref{induced-cycles} (ii)
implies
$$
Z_{\line G} = Z^{local}_{\line G} + Z^{global}_{\line G}.
$$

To simplify the presentation \eqref{transpose-cokernel-presentation},
note that for a vertex $vw$ of $\line G$, the bond
$$
b_{\line G}(vw) = f^t(vw,v) + f^t(vw,w)
$$
lies in $\im(f^t)$, and consequently, $B_{\line G} \subset \im(f^t)$.
Note also that the decomposition \eqref{disjoint-edge-decomposition}
leads to a family of compatible direct sum decompositions
\begin{align*}
  \ZZ^{E_{\line G}} &= \bigoplus_{v \in V_G} \ZZ^{E(K_d^{(v)})} \\
  Z^{local}_{\line G} &= \bigoplus_{v \in V_G} Z_{K_d^{(v)}} \\
  \im(f^t) & = \bigoplus_{v \in V_G} B_{K_d^{(v)}}.
\end{align*}
This gives the simplified presentation
\begin{equation}
\label{simplified-presentation}
\begin{aligned}
\coker(f^t) &= \left( \bigoplus_{v \in V_G} \ZZ^{E(K^{(v)}_d)} / 
               \left( B_{K_d^{(v)}} +  Z_{K_d^{(v)}} \right) \right)
               / Z^{global}_{\line G} \\
&= \left( \bigoplus_{v \in V_G} K(K^{(v)}_d) \right) / Z^{global}_{\line G}.
\end{aligned}
\end{equation}

We use this presentation to prove the following lemma, which together
with Proposition~\ref{cokernel-description} immediately implies
Theorem~\ref{regular-result}.
\begin{lemma}
\label{kernel-description}
For a connected $d$-regular graph $G$,
$$
\ker(f) \cong \ZZ_d^{\beta(G)-2} \oplus C
$$
where $C:=\coker(f)$ is as described in Theorem~\ref{regular-result}.
\end{lemma}
\begin{proof}
  We claim that it suffices to prove these two bounds on $\ker(f)$:
  \begin{itemize}
  \item[](i) There is a surjection $\ker(f) \twoheadrightarrow
    \ZZ_d^{\beta(G)-2}$ and,
  \item[](ii) $\ker(f)$ can be generated by $\beta(G)-1$ elements.
  \end{itemize}

  To see this claim, note that since $\ker(f)$ is all $d$-torsion by
  Corollary~\ref{f-torsion-corollary}, assertion (ii) would imply a
  surjection $\ZZ_d^{\beta(G)-1} \twoheadrightarrow \ker(f)$.
  Together with (i), this would imply $\ker(f) \cong
  \ZZ_d^{\beta(G)-2}\oplus C'$ for some cyclic group $C'$.  But then
  exactness of the sequence
$$
0 \rightarrow 
\underbrace{\ker(f)}_{\ZZ_d^{\beta(G)-2} \oplus C'} \rightarrow 
K(\line G) \overset{f}{\rightarrow} 
K(\sd G) \rightarrow 
\underbrace{\coker(f)}_{C} \rightarrow 
0
$$
forces
$$
\left( d^{\beta(G)-2} |C'| \right) |K(\sd G)| = |K(\line G)| |C|.
$$ 
From this equation and equation \eqref{Sachs'-theorem} one deduces $|C'|=|C|$.
Since both $C'$ and $C$ are cyclic, this means $C' \cong C$, as desired.

In the proofs of assertions (i) and (ii), one uses the fact that $\ker(f)=\coker(f^t)$.
Moreover, setting $n:=|V_G|$,  one can rewrite the direct sum from \eqref{simplified-presentation} as
\begin{equation}
\label{local-direct-sum}
\bigoplus_{v \in V_G} K(K^{(v)}_d) 
  \cong   \bigoplus_{v \in V_G} \ZZ_d^{d-2} 
  \cong  \ZZ_d^{n(d-2)}.
\end{equation}

For assertion (i), we use some easy numerology.  Note that
$Z^{global}_{\line G}$ can be generated by $\beta(G)$ elements,
and also that 
$$
\beta(G) = |E_G| - |V_G| + 1 = \frac{dn}{2} - n +1 =  \frac{n(d-2)}{2} + 1
$$
so that 
$$
n(d-2)-\beta(G) = \beta(G)-2.
$$
Since it is easily seen that that any quotient of an abelian group $\ZZ_d^a$ by a subgroup that can be generated by $b$ elements
will have a surjection to $\ZZ_d^{a-b}$, one can apply this with $a=n(d-2)$ and $b=\beta(G)$ to
the presentation \eqref{simplified-presentation}, and conclude that there is a
surjection $\coker(f^t) \twoheadrightarrow \ZZ_d^{\beta(G)-2}$.

For assertion (ii), the idea will be to start with the 
$$
n(d-2)=2(\beta(G)-1)
$$
generators in \eqref{local-direct-sum}, and use (all but one of) the $\beta(G)$
generating global cycles 
in $Z^{global}_{\line G}$ to rewrite them in terms of other generators, with $\beta(G)-1$ 
generators left.  This will be
achieved by removing the vertices from $G$ one at a time in a certain order, in order to control
the rewriting process.

To this end, order the vertices $V_G$ as $v_1,v_2,\ldots,v_n$ in such a way that the
vertex-induced subgraphs 
$$
\begin{aligned}
G_i &:=G \setminus \{v_1, v_2,\ldots,v_{i-1}\}  \\
    & (\text{so }G_1:=G, \text{ and }G_n \text{ has one vertex }v_n)
\end{aligned}
$$
satisfy
$$
d_i:=\deg_{G_i}(v_i) < d \text{ for } i \geq 2.
$$
For each $i \geq 1$, partition the $d_i$ neighbors $v_i$ in $G_i$
into blocks $A_1,A_2,\ldots,A_{c_i}$ according to the connected components
of $G_{i+1}$ in which they lie.  The number $c_i$ of such components coincides with the
number of connected components in $G_{i+1}$ into which the connected component
of $v_i$ in $G_i$ splits after removing $v_i$.
Define
$$
\Delta_i:= d_i - c_i = \beta(G_i) - \beta(G_{i+1}),
$$
where the last equality follows from the Euler relation for graphs:
$$
|V_G|-|E_G| = |\{\text{connected components of }G\}| - |\beta(G)|.
$$
Consequently,
$$
\Delta_1 + \Delta_2 + \cdots + \Delta_{n-1} = \beta(G_1)-\beta(G_n) = \beta(G).
$$
Our goal will then be to find $\Delta_i$ minimal generators of
\eqref{local-direct-sum} to remove at
each stage $i \geq 2$ (and at the first stage $i=1$, remove one fewer, that is,
$\Delta_1-1=d-2$ of them).
This would leave a generating set for $\coker(f^t)$ of cardinality
$n(d-2)-(\beta(G)-1) = \beta(G)-1$, as desired.

For $i \geq 2$, inside the clique $K^{(v_i)}_d$ local to $v_i$, choose
a forest $F_i$ of
edges having $c_i$ components which are spanning trees on each of the subsets
$\{v_i x: x \in A_j\}$ for $j=1,2,\ldots,c_i$.  Note that 
$$
|F_i|= \sum_{j=1}^{c_i} (|A_j|-1) = d_i - c_i = \Delta_i.
$$ 
Also note that the forest $F_i$ manages to avoid touching at least one
vertex in the $d$-clique $K^{(v_i)}_d$, namely any vertex of the form
$v_i v_k$ in which $\{ v_i v_k \} \in E_G$ and $k < i$;  there will exist
at least one such $k$ since by construction, 
$\deg_{G_i}(v_i) =d_i < d = \deg_{G}(v_i)$.

Hence by Proposition~\ref{complete-graph}, the edges in $F_i$ give $\Delta_i$
generators that could be completed to a set of $d-2$ minimal generators for
$K(K^{(v_i)}_d) \cong \ZZ_d^{d-2}$.  Each of these generators in $F_i$
can be re-written, using a cycle in $Z^{global}_{\line G}$, in terms of
generators from $K(K^{(v_k)})$'s that have $k > i$, as follows.  Given
any edge $(v_i x,v_i x')$ in $F_i$, there is a path from $x$ to $x'$
in $G_{i+1}$ (because $x, x'$ lie in the same component of $G_{i+1}$ by
construction), and hence a directed cycle $C$ in $G_i$ going
from $v_i$ to $x$ then through this path to $x'$ and back to $v_i$.
The global cycle $z(\line C)$ allows one to rewrite $(v_i x,v_i x')$ as desired.

The only difference for $i=1$ is that, even when $\Delta_1=d-1$ 
(that is, when $c_1=1$), in this situation choose $F_i$ to have at most
$d-2$ edges (that is, remove any edge from the forest $F_1$ if $c_1=1$).  This modification ensures
that one can still apply Proposition~\ref{complete-graph} and rewrite all
of the generators of $K(K_d^{(v_1)})$ corresponding to the edges of $F_1$.
\end{proof}

\begin{remark}
  One should remark that for a connected, $d$-regular graph $G$, the
  extra hypothesis in Theorem~\ref{regular-result} that $G$ is
  $2$-edge-connected is well-known (see, e.g., \cite{Hartsfield}) to be
  superfluous when $d$ is even: a connected graph $G$ with all
  vertices of even degree cannot have a cut-edge, as the two
  components created by the removal of this edge would each be graphs
  having exactly one vertex of odd degree, an impossibility.

  However, when $d$ is odd, the extra hypothesis of
  $2$-edge-connectivity need not follow.  For example, the $3$-regular
  graph shown in Figure~\ref{3-regular} is connected, but not
  $2$-edge-connected.

\begin{figure}
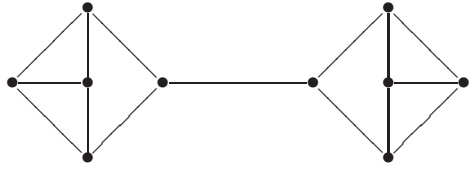

\[
\xygraph{
!{<0cm,0cm>;<0cm,1cm>:<1cm,0cm>::}
!{(0,-3) }*{\bullet}="a"
!{(0,-2) }*{\bullet}="b"
!{(0,-1) }*{\bullet}="c"
!{(0,1) }*{\bullet}="d"
!{(0,2) }*{\bullet}="e"
!{(0,3) }*{\bullet}="f"
!{(-1,2) }*{\bullet}="g"
!{(-1,-2) }*{\bullet}="h"
!{(1,2) }*{\bullet}="i"
!{(1,-2) }*{\bullet}="j"
"a"-"b"-"h"-"a"-"j"-"b"
"h"-"c"-"j"
"c"-"d"
"g"-"d"-"i"
"g"-"e"-"i"
"g"-"f"-"i"
"f"-"e"
}
\]
\caption{A 3-regular connected graph which is not 2-edge-connected.}
\label{3-regular}
\end{figure}

\end{remark}

\section{Proof of Corollary~\ref{regular-nonbipartite-corollary}}
\label{regular-nonbipartite-section}

In this section we prove
Corollary~\ref{regular-nonbipartite-corollary}. Informally, the
corollary states that critical group of $G$ determines the critical
group of $\line G$ in a simple way.

\vskip.1in
\noindent
{\bf Corollary~\ref{regular-nonbipartite-corollary}.}
{\it
For $G$ a simple, connected, $d$-regular graph with $d \geq 3$ which is nonbipartite,
after expressing uniquely  
$$
K(G) \cong \bigoplus_{i=1}^{\beta(G)} \ZZ_{d_i}
$$ 
with $d_i$ dividing $d_{i+1}$, one has
\begin{equation}
\label{algebraic-Sachs}
K(\line G) \cong \left( \bigoplus_{i=1}^{\beta(G)-2} \ZZ_{2dd_i} \right)
\oplus  \begin{cases}
\ZZ_{2d_{\beta(G)-1}} \oplus \ZZ_{2d_{\beta(G)}} & \text{ for }|V|\text{ even,}\\
\ZZ_{4d_{\beta(G)-1}} \oplus \ZZ_{d_{\beta(G)}} & \text{ for }|V|\text{ odd.}
\end{cases}
\end{equation}
}

\begin{proof}
  Let $K:=K(\line G)$, and fix a prime $p$.  Our goal is to show that
  the $p$-primary component of $K$ matches that of the group on the
  right side of \eqref{algebraic-Sachs}.

  The hypotheses of the theorem allow one to apply the nonbipartite
  cases of Theorem~\ref{divisibility-theorem} and
  Theorem~\ref{regular-result}.  The former asserts that
\begin{equation}
  \label{nonbipartite-divisibility}
  K/p^{k(p)}K \cong \ZZ_{p^{k(p)}}^{\beta(G)-2} \oplus
  \begin{cases}
    0 & \text{ for }p\text{ odd,}\\
    \ZZ_2^2 &\text{ for }p=2, |V|\text{ even,}\\
    \ZZ_4 &\text{ for }p=2, |V|\text{ odd,}
  \end{cases}
\end{equation}
while the latter gives an exact sequence
\begin{equation}
\label{nonbipartite-exact-sequence}
0 \rightarrow \ZZ_d^{\beta(G)-2} \oplus \ZZ_{\gcd(2,d)}
  \rightarrow K
  \rightarrow K(\sd G)
  \rightarrow \ZZ_{\gcd(2,d)} \rightarrow 0.
\end{equation}

In analyzing the $p$-primary component $\Syl_p(K)$, it is convenient to define
the {\it type} of a finite abelian $p$-group $A$ as
the unique integer partition $\nu=(\nu_1 \geq \nu_2 \geq \cdots)$
for which 
$A \cong \bigoplus_{i \geq 1} \ZZ_{p^{\nu_i}}.$
Let $\mu, \lambda$ denote the types of $\Syl_p (K(G)), \, \Syl_p(K)$,
where we think of both $\mu, \lambda$ as partitions with $\beta(G)$ parts
(allowing some parts to be $0$).  Note that
Proposition~\ref{Lorenzini's-theorem} asserts, in this language,
that $\Syl_p(K(\sd G))$ has type $\mu$ for $p$ odd and
type $\mu+(1^{\beta(G)})$ for $p=2$.

A basic fact from the theory of Hall polynomials \cite[Chapter II
Section 9]{Macdonald} says that there exist short exact sequences of
abelian $p$-groups 
$$
0 \rightarrow A \rightarrow B \rightarrow C
\rightarrow 0
$$ 
in which $A,B,C$ have types $\nu,\lambda,\mu$, respectively, if and
only if the {\it Littlewood-Richardson (or LR) coefficient}
$c_{\mu,\nu}^\lambda$ does not vanish.  The combinatorial rephrasing
of this {\it LR-condition} is as follows: There must exist at least
one {\it column-strict tableau} (which we will call an {\it LR
  tableau}) of the skew-shape $\lambda/\mu$ having content $\nu$, for
which the word obtained by reading the tableau (in English notation)
from right-to-left in each row, starting with the top row, is {\it
  Yamanouchi}.  Here the Yamanouchi condition means that within each
initial segment of the word, and for each value $i \geq 1$, the number
of occurrences of $i+1$ is at most the number of occurrences of $i$.
See \cite[Chapter I \S 9]{Macdonald} and \cite[Appendix \S
A1.3]{Stanley} for more on these notions.

Suppose that $p$ is odd. Then $k(p)$ is the largest power such that
$p^{k(p)}$ divides $d$, so taking the $p$-primary components in
\eqref{nonbipartite-exact-sequence}, we obtain the following short exact
sequence:
\begin{equation}
0 \longrightarrow \underbrace{\ZZ_{p^{k(p)}}^{\beta(G)-2}}_{\text{ type }\nu=(k(p)^{\beta(G)-2})}
  \longrightarrow \underbrace{\Syl_p K}_{\text{ type }\lambda}
  \longrightarrow \underbrace{\Syl_p( K(\sd G) )}_{\text{ type }\mu}
  \longrightarrow 0
\end{equation}
where $\lambda$ has at most $\beta(G)-2$ nonzero parts by
\eqref{nonbipartite-divisibility}.  Since nonvanishing of the
LR-coefficient $c^{\lambda}_{\mu,\nu}$ forces $\mu \subset \lambda$,
it must be that $\mu$ also has at most $\beta(G)-2$ nonzero parts.
Furthermore, one can check that column-strictness together with the
Yamanouchi condition on the reading word of an LR-tableau of shape
$\lambda/\mu$ and content $\nu=(k(p)^{\beta(G)-2})$ uniquely determine
the tableau: It must have each entry in row $i$ equal to $i$ for
$i=1,2,\ldots,\beta(G)-2$.  This forces $\lambda_i=\mu_i+k(p)$ for
$i=1,2,\ldots,\beta(G)-2$, and hence $\lambda$ agrees with the type of
the $p$-primary component on the right side of
\eqref{algebraic-Sachs}.

Suppose that $p=2$, so that $2^{k(p)-1}$ divides $d$, but $2^{k(p)}$
does not.

When $d$ is odd we have that $k(p)=1$. On the other hand, taking the
$2$-primary components in \eqref{nonbipartite-exact-sequence} shows
that $\Syl_2 K \cong \Syl_2(K(\sd G))$, so
$\lambda=\mu+(1^{\beta(G)})$.  Since $d$ is odd, $|V|$ must be even
(as the $d$-regularity of $G$ forces $d|V|=2|E|$), so this $\lambda$
again agrees with the type of the $2$-primary component on the right
side of \eqref{algebraic-Sachs}.

If $d$ is even, the $2$-primary components in
\eqref{nonbipartite-exact-sequence} form the following exact sequence
\begin{equation}
0 \longrightarrow \underbrace{\ZZ_{2^{k(p)-1}}^{\beta(G)-2} \oplus \ZZ_2}_{\text{ type }\nu=((k(p)-1)^{\beta(G)-2},1)}
  \longrightarrow \underbrace{\Syl_2 K}_{\text{ type }\lambda}
  \longrightarrow \underbrace{\Syl_2( K(\sd G) )}_{\text{ type }\mu+(1^{\beta(G)})}
  \overset{\pi}{\rightarrow} \ZZ_2 
  \longrightarrow 0
.
\end{equation}
This can be truncated to the following {\it short} exact sequence involving $\ker \pi$:
\begin{equation}
0 \longrightarrow \underbrace{\ZZ_{2^{k(p)-1}}^{\beta(G)-2} \oplus \ZZ_2}_{\text{ type }\nu=(k(p)-1)^{\beta(G)-2},1)}
  \longrightarrow \underbrace{\Syl_2 K}_{\text{ type }\lambda}
  \longrightarrow \underbrace{\ker\pi}_{\text{ type }\hat\mu}
  \longrightarrow 0 
\end{equation}
for some partition $\hat\mu$, and where the last two parts
$
(\lambda_{\beta(G)-1},\lambda_{\beta(G)})
$ 
in $\lambda$ are
either $(1,1)$ or $(2,0)$ by \eqref{nonbipartite-divisibility},
depending on the parity of $|V|$.  

The short exact sequence 
$$
0 \rightarrow \ker\pi \rightarrow \Syl_2(
K(\sd G) ) \overset{\pi}{\rightarrow} \Z_2 \rightarrow 0
$$ 
shows that $\hat\mu$ is obtained from $\mu+(1^{\beta(G)})$ by removing
one square; we claim that $\hat\mu$ can have at most $\beta(G)-1$
nonzero parts, and hence this square must be removed from the {\it
  last} row, that is, $\hat\mu=\mu+(1^{\beta(G)-1},0)$.  The reason
for this claim is that, since the LR-coefficient
$c^{\lambda}_{\hat\mu,\nu} \neq 0$, the LR-condition forces
$$
\sum_{i \geq \beta(G)-1} \lambda_i \geq 
\sum_{i \geq \beta(G)-1} \hat\mu_i +  \sum_{i \geq \beta(G)-1} \nu_i. 
$$
As $\sum_{i \geq \beta(G)-1} \lambda_i=2$ in both cases for the
parity of $|V|$, and $\sum_{i \geq \beta(G)-1} \nu_i =1$,
this forces $\sum_{i \geq \beta(G)-1} \hat\mu_i \leq 1$.  This
implies $\hat\mu$ can have at most $\beta(G)-1$ nonzero parts, as claimed.

Once one knows $\hat\mu$ takes this form, and since
\eqref{nonbipartite-divisibility} fixes the shape of $\lambda$
in its last two rows $\beta(G)-1, \beta(G)$, any 
LR-tableau of shape $\lambda/\mu$ and content $\nu=(k(p)-1)^{\beta(G)-2},1)$
is completely determined by column-strictness and the Yamanouchi condition:
It must have its unique entry equal to $\beta(G)-1$ lying in the unique of cell
of $\lambda/\mu$ within the last two rows, while all of its entries in row $i$
are all equal to $i$ for $i=1,2,\ldots,\beta(G)-2$.
This again forces $\lambda_i=\mu_i+k(p)$ for $i=1,2,\ldots,\beta(G)-2$,
and means that $\lambda$ again matches the type
of the $2$-primary component on the right side of \eqref{algebraic-Sachs}.
\end{proof}

\begin{remark}
  In light of what Corollary~\ref{regular-nonbipartite-corollary} says
  about $K:=K(\line G)$ for {\it nonbipartite} regular graphs, one
  might wonder what can be deduced for {\it bipartite} regular graphs
  using Theorems~\ref{divisibility-theorem} and
  \ref{regular-result}. We discuss this briefly here.

  Fixing a prime $p$, define $k$ to be the largest exponent such that
  $p^k$ divides $d$, and let $\Syl_p(K(G))$ have type $\mu$.  Then the
  $p$-primary components in the bipartite case of
  Theorem~\ref{regular-result} form the following exact sequence:
\begin{equation}
\label{bipartite-regular-sequence}
0 \rightarrow \underbrace{\ZZ_{p^{k}}^{\beta(G)-1}}_{\text{type }\nu=(k^{\beta(G)-1})}
  \longrightarrow \underbrace{ \Syl_p K}_{\text{type }\lambda}
  \longrightarrow \underbrace{\Syl_p( K(\sd G) )}_{
    \tiny\begin{cases}
      \text{type }\mu &\text{ if }p \neq 2\\
      \text{type }\mu+(1^{\beta(G)}) &\text{ if }p=2
    \end{cases}
}
  \overset{\pi}{\longrightarrow} \underbrace{\ZZ_{p^k}}_{\text{type }(k)}
  \rightarrow 0
\end{equation}

As a consequence, $\Syl_p(K)$ will be uniquely determined by
$\Syl_p(K(G))$ whenever $p$ does not divide $d$, since
then $k=0$ and \eqref{bipartite-regular-sequence} shows
$\Syl_p(K) \cong \Syl_p(K(\sd G))$ in this case.
However, in general, the structures of $\ker(\pi)$ and of $\Syl_p(K)$ seem less clear.
Even using the extra information from Theorem~\ref{divisibility-theorem}
that $K/p^{k(p)}K \cong \ZZ_{p^{k(p}}^{\beta(G)-1} \oplus \ZZ_{\gcd(p^{k},|V|)}$,
where $k(p)$ is the largest power such that $p^{k(p)}$ divides $2d$,
along with the LR-rule, the structures of the various terms in the
sequence are not uniquely determined.

\begin{question}
\label{bipartite-regular-question}
When $G$ is a simple, bipartite, regular graph, what
more can be said about the structure of
$K:=K(\line G)$ in relation to that of $K(G)$?
\end{question}

\end{remark}

\section{Proof of Theorem~\ref{semiregular-result}}
\label{semi-regular}

Let $G=(V,E)$ be a semiregular bipartite graph with vertex bipartition
$V = V_1 \sqcup V_2$, such that vertices in $V_i$ have degree $d_i$.
In this section we prove our analogue of 
Theorem~\ref{regular-result} for semiregular graphs. Recall that this
is motivated by Cvetkovi\'c's formula
\eqref{Cvetkovic's-theorem} for the spanning tree number of $\line G$:
$$
\kappa(\line G) = \frac{(d_1+d_2)^{\beta(G)}}{d_1 d_2}  
\left(\frac{d_1}{d_2}\right)^{|V_2|-|V_1|} \kappa(G). 
$$
We recall here the statement of Theorem~\ref{semiregular-result}.
\begin{semiregtheorem*}
Let $G$ be a connected bipartite $(d_1,d_2)$-semiregular graph $G$.
Then there is a group homomorphism 
$$
K(\line G) \overset{g}{\rightarrow} K(G)
$$
whose kernel-cokernel exact sequence
\begin{equation}
0 
\rightarrow \ker(g)
\rightarrow K(\line G)
\overset{g}{\rightarrow} K(G)
\rightarrow \coker(g)
\rightarrow 
0
\end{equation}
has
\begin{itemize}
\item
$\coker(g)$ all $\lcm(d_1,d_2)$-torsion, and
\item
$\ker(g)$ all $\frac{d_1+d_2}{\gcd(d_1,d_2)} \lcm(d_1,d_2)$-torsion.
\end{itemize}
\end{semiregtheorem*}
The proof of this result is analogous to that of
Theorem~\ref{regular-result}; for this reason, some proofs here are
either abbreviated or only sketched.  Note also that this theorem is
less precise than Theorem~\ref{regular-result}, partly out of
necessity: Examples~\ref{regular-complete-bipartite-example} and
\ref{complete-bipartite-example} below show that the morphism $g:
K(\line G)\rightarrow K(G)$ appearing in the theorem is nearly
surjective in some cases, but is the zero morphism in some other
cases!

\subsection{Defining the morphism $g$}
\label{defining-g-subsection}

We define $g$ similarly to the map $f$ from
Definition~\ref{f-definition}. Let
$$
\begin{aligned}
\lambda &:=\lcm(d_1,d_2) \\ 
\gamma &:=\gcd(d_1,d_2).
\end{aligned}
$$
As a notational convenience, denote typical vertices in $V_1$
(respectively, $V_2$) by $a$'s (respectively, $b$'s) with subscripts
or primes.

\begin{definition}
  For a semiregular bipartite graph $G$, let $g: \RR^{E_{\line G}}
  \rightarrow \RR^{E_G}$ be defined $\RR$-linearly by
  $$
  \begin{aligned}
    g(ab,ba') & =\frac{\lambda}{d_2} \left( (a,b) +  (b,a') \right) \\
    g(ba,ab') & =\frac{\lambda}{d_1} \left( (b,a) +  (a,b') \right).
  \end{aligned}
  $$
  Equivalently, the adjoint map $g^t$ is defined by
  $$
  g^t (a,b)=\frac{\lambda}{d_1}\sum_{b_i \in N(a)}(b_i a,ab) 
  + \frac{\lambda}{d_2}\sum_{a_j \in N(b)}(ab,b a_j),
  $$
  where $N(v)$ denotes the set of vertices adjacent to $v$ in $G$.
\end{definition}

\begin{remark}
\label{regular-bipartite-case}
In the special case when $G$ is not only semiregular bipartite, but
actually regular, so $d_1=d_2 =\lambda =\gamma$, one can easily check
that the map $g$ coincides with the composite map $h \circ f$
$$
\RR^{E_{\line G}} \overset{f}{\longrightarrow}\RR^{E_{\sd G}}
\overset{h}{\longrightarrow}\RR^{E_G}
$$
where $f$ is the map from Theorem~\ref{regular-result} defined in
Definition~\ref{f-definition}, and $h$ was defined in
Example~\ref{subdivision-example}.
\end{remark}

\begin{proposition}
  If $G$ is a semiregular bipartite graph, then $g: \ZZ^{E_{\line G}}
  \rightarrow \ZZ^{E_G}$ is a morphism of the associated rational
  orthgonal decompositions, and hence induces a group homomorphism $g:
  K(\line G) \rightarrow K(G)$.
\end{proposition}

\begin{proof} 
  By Lemma \ref{induced-cycles} (ii), it is enough to show that $g$
  takes global and local cycles in $Z_{\line G}$ to cycles in $Z_G$,
  and that $g^t$ takes cycles in $Z_G$ to cycles in $Z_{\line G}$.

  First, one can check that $g$ maps all local cyles to $0$.  Each
  global cycle is by definition of the form $z(\line C)$ where $C$ is
  a directed cycle of $G$, and one checks that
  $$
  g(z(\line C)) =  \left( \frac{\lambda}{d_1} + \frac{\lambda}{d_2} \right) z(C).
  $$

  On the other hand, one checks that $g^t(z(C))$ can be rewritten as a
  sum of $\lambda$ cycles $\zeta_i$ in $Z_{\line G}$, each $\zeta_i$
  being twice the length of $C$, and in which every other vertex on
  $\zeta_i$ corresponds to an edge occurring in $C$.
\end{proof}

\subsection{Analyzing its kernel and cokernel}
\label{semiregular-kernel-cokernel-subsection}

\begin{proposition}
The map 
$$
g^t g: K(\line G) \rightarrow K(\line G)
$$
coincides with scalar multiplication by $\frac{d_1+d_2}{\gamma}
\lambda$.  Consequently, $\ker(g)$ is $\frac{d_1+d_2}{\gamma}
\lambda$-torsion.
\label{g-scalarmult}
\end{proposition}

\begin{proof}
  For any edge $ab,ba'$ in $E_{\line G}$, use the definitions of $g$
  and $g^t$ to write
\begin{equation}
\label{eqngtg}
\begin{aligned}
  g^tg(ab,ba')  =&   \frac{\lambda^2}{d_1 d_2} \sum_{b_i \in N(a)} (b_ia,ab)  
  + \frac{\lambda^2}{d_2^2} \sum_{a_j \in N(b)} (ab,ba_j)   \\
  & + \frac{\lambda^2}{d_2^2} \sum_{a_j \in N(b)} (a_jb,ba') 
  + \frac{\lambda^2}{d_1 d_2} \sum_{b_k \in N(a')} (ba',a'b_k).
\end{aligned}
\end{equation}
For the first and fourth term, one has
$$
\begin{aligned}
\sum_{b_i \in N(a)} (b_ia,ab ) &= \sum_{a_j \in N(b)} (ab,ba_j) \mbox{ mod } B_{\line G} \\
\sum_{b_k \in N(a')} (ba',a'b_k) & = \sum_{a_j \in N(b)} (a_jb,ba') \mbox{ mod } B_{\line G}.
\end{aligned}
$$
Substituting these expressions into equation \eqref{eqngtg}, grouping
like terms, and using the identity \newline $d_1 d_2 = \lambda \gamma$
gives
$$
g^tg(ab,ba') = \frac{d_1 + d_2}{\gamma} \cdot \frac{\lambda}{d_2} 
\left(\sum_{a_j \in N(b)} (ab,ba_j) + \sum_{a_j \in N(b)} (a_jb,ba') \right) 
\mbox{ mod } B_{\line G},
$$
which then can be rewritten, using the $d_2$ triangular cycles 
$$
(ab,ba_j) + (a_jb,ba') + (a'b,ba) \in Z_{\line G},
$$
as
\begin{align*}
  g^tg(ab,ba') & = \frac{(d_1 + d_2)}{\gamma} \cdot
  \frac{\lambda}{d_2} \left( d_2 (ab,ba') \right)
  \mbox{ mod } B_{\line G} + Z_{\line G} \\
  & = \frac{(d_1 + d_2)}{\gamma} \lambda (ab,ba') \mbox{ mod }
  B_{\line G} + Z_{\line G}.\qedhere
\end{align*}
\end{proof}

\begin{remark}
As in Proposition~\ref{f-scaling-prop}, one can show that the other
map $g g^t: K(G) \rightarrow K(G)$ also coincides with the scalar multiplication
by $\frac{d_1+d_2}{\gamma} \lambda$, and hence that $\coker(g)$ is also
$\frac{d_1+d_2}{\gamma} \lambda$-torsion.  However, we omit this proof,
since we are about to show the {\it stronger} assertion that $\coker(g)$ is $\lambda$-torsion.
\end{remark}

\vskip.2in
\noindent
{\it Proof of Theorem~\ref{semiregular-result}.}
In light of Proposition~\ref{g-scalarmult}, it only remains to show that
$\coker(g)$ is $\lambda$-torsion. Given any edge $ab \in E_G$, one has
$$
\begin{aligned}
\lambda (a,b) & = \lambda (a,b) + 
                 \frac{\lambda}{d_2} \sum_{a_j \in N(b)} (b,a_j) \mbox{ mod } B_{G} \\
              & = \frac{\lambda}{d_2} \sum_{a_j \in N(b)} (a,b) + (b,a_j) \mbox{ mod } B_{G} \\
              & = g \left( \sum_{a_j \in N(b)} (ab,ba_j) \right) \mbox{ mod } B_{G} .
\end{aligned}
$$
Consequently $\lambda (a,b)$ lies in $\im(g)+B_{G}$, so it
is zero in $\coker(g) : = \ZZ^{E_G}/ \left( \im(g) + B_G + Z_G \right)$.
$\qed$

\vskip .2in

Unlike the map $f$ from Section~\ref{regular-result-section}, it is
hard to be much more precise about the exact nature of cokernel and
kernel of $g$.  The following two families of examples demonstrate two
extremes of behavior for how tightly or loosely the map $g$ ties
together $K(\line G)$ and $K(G)$ for semiregular bipartite graphs $G$.

\begin{example}
\label{regular-complete-bipartite-example}
Assume $G$ is not only bipartite semiregular, but actually $d$-regular
(i.e., $d_1=d_2=d$).  Then $g: K(\line G) \rightarrow K(G)$ is nearly
surjective, in the sense that $\coker(g)$ is a quotient of $\ZZ_d$.
To see this, recall from Remark~\ref{regular-bipartite-case} that in
this case, $g = h \circ f$ where $h, f$ were defined in
Example~\ref{subdivision-example} and Definition~\ref{f-definition}.
Since $h: K(\sd G) \rightarrow K(G)$ is surjective, it induces a
surjection
$$
\coker(f) := K(\sd G)/\im(f) \longrightarrow K(G)/\im(h \circ f) =: \coker(g).
$$
But Theorem~\ref{regular-result} says that $\coker(f) = \ZZ_d$ in this situation.
\end{example}

\begin{example}
\label{complete-bipartite-example}
For the complete bipartite graph $G=K_{n_1,n_2}$, the structures of
the critical groups of $G$ and $\line G$ have been determined through
manipulations of their Laplacian matrices (see Lorenzini
\cite{Lorenzini} and Berget \cite{REU}, respectively):

\begin{equation}
\begin{aligned}
K(K_{n_1,n_2}) & \cong 
\ZZ_{n_1}^{n_2-2} \oplus \ZZ_{n_2}^{n_1-2} \oplus \ZZ_{n_1n_2}, \\
K(\line K_{n_1,n_2}) & \cong 
   \ZZ_{n_1(n_1+n_2)}^{n_1-2} \oplus \ZZ_{n_2(n_1+n_2)}^{n_2-2} 
         \oplus \ZZ_{n_1+n_2}^{(n_1-2)(n_2-2)+1}.
\end{aligned}
\nonumber
\end{equation}


In principle, the structures of these groups allow nonzero
homomorphisms between them for all values of $n_1, n_2$.  However, we
claim that whenever $\gamma=\mbox{gcd}(n_1,n_2)=1$, the map $K(\line
G) \overset{g}{\rightarrow} K(G)$ will be the zero morphism. In this
case, $\ker(g) = K(\line G)$ and $\coker(g) = K(G)$.

To see this claim, let $(ab,ba')$ be a fixed edge in $E_{\line G}$.
Note that
$$
\begin{aligned}
\lambda  &= d_1 d_2 = n_1 n_2, \\
d_1 &= n_2, \,\, d_2 =n_1.
\end{aligned}
$$

Then for each $b_j \in V_2$, one has 
$$
\frac{1}{d_1}\left( g(ab,ba')+g(a'b_j,b_ja) \right)
=
ab+ba'+a'b_j+b_ja \in Z_G.
$$
On the other hand,
$$
\begin{aligned}
\sum_{b_{j} \in V_2} \frac{1}{d_1}\left( g(ab,ba')+g(a'b_j,b_ja) \right)
&=
g(ab,ba') + \sum_{b_{j} \in V_2}\frac{1}{d_1}\frac{\lambda}{d_2}(a'b_j+b_ja) \\
&= 
g(ab,ba') +\sum_{b_{j} \in V_2}a'b_j+\sum_{b_{j} \in V_2}b_ja  \\
&=
g(ab,ba') \mbox{ mod } B_G.
\end{aligned}
$$
Combining these two statements gives us $g(ab,ba') = 0$ mod
$Z_G + B_G$.  By symmetry, one also has $g(ba,ab') = 0$ mod $Z_G + B_G$.  
It follows that $g$ is the zero morphism.
\end{example}

\begin{remark}
Note that Theorem~\ref{divisibility-theorem} provides
convenient information about $\Syl_p (K)$ for $K:= K(\line G)$
when $G$ is $(d_1,d_2)$-semiregular: If $k(p)$ denotes
the largest power $p^{k(p)}$ dividing $d_1+d_2$,
then 
\begin{equation}
\label{convenient-semiregular-info}
K/p^{k(p)}K \cong \ZZ_{p^{k(p)}}^{\beta(G)-1} \oplus \ZZ_{\gcd(p^{k(p)},|V|)}.
\end{equation}
However, even in conjunction with Theorem~\ref{semiregular-result},
this does not appear to determine the structure of $K(\line G)$
uniquely in terms of the structure of $K(G)$.  Thus we are led to the
following generalization of Question~\ref{bipartite-regular-question}:

\begin{question}
\label{bipartite-semiregular-question}
When $G$ is a simple, semiregular bipartite graph, what more can be
said about the structure of $K:=K(\line G)$ in relation to that of
$K(G)$?
\end{question}

\end{remark}

\section{Examples}\label{example-section}

\subsection{The complete graph $K_n$}
\label{complete-graph-section}
Proposition~\ref{complete-graph} can be rephrased as asserting that
$$
K(K_n) \cong \ZZ_n^{n-2} \oplus \ZZ_1^{\beta-n+2}
$$
where $\beta:=\beta(K_n)=\binom{n-1}{2}$.  Since $K_n$ is nonbipartite
for $n \geq 3$, and contains an even length cycle for $n \geq 4$,
Corollary~\ref{regular-nonbipartite-corollary} immediately implies the
following:
\begin{corollary}
\label{complete-graph-corollary}
For $n \geq 4$, the line graph $\line K_n$ of the critical group of
the complete graph $K_n$ has the form
\begin{align*}\label{complete-graph-equation}
  K(\line K_4) &= \ZZ_{24} \oplus \ZZ_8 \oplus \ZZ_2 \\
  K(\line K_n) &= \ZZ_{2(n-1)n}^{n-2} \oplus \ZZ_{2(n-1)}^{\beta-n} \oplus 
  \begin{cases}
    \ZZ_2^2 &\text{ for even }n > 5,\\
    \ZZ_4 & \text{ for odd }n \geq 5.
  \end{cases}
\end{align*}

\end{corollary}
\subsection{The complete bipartite graph $K_{n_1,n_2}$ }
\label{complete-bipartite-graph-section}

As mentioned in Example~\ref{complete-bipartite-example}, the critical
groups of the complete bipartite graph $K_{n_1,n_2}$ and its line
graph $\line K_{n_1,n_2}$ have the following forms:
$$
\begin{aligned}
K(K_{n_1,n_2})  &\cong 
  \ZZ_{n_1}^{n_2-2} \oplus \ZZ_{n_2}^{n_1-2} \oplus \ZZ_{n_1n_2}\\
K(\line K_{n_1,n_2}) &\cong
   \ZZ_{n_1(n_1+n_2)}^{n_1-2} \oplus \ZZ_{n_2(n_1+n_2)}^{n_2-2} 
         \oplus \ZZ_{n_1+n_2}^{(n_1-2)(n_2-2)+1}
\end{aligned}
$$
(see Lorenzini \cite{Lorenzini} and Berget \cite{REU}, respectively).
 
In addition, Example~\ref{complete-bipartite-example} showed that the
map $g$ in the exact sequence in Theorem~\ref{semiregular-result} is
sometimes the zero morphism and hence is not always useful for
determining the structure of $K(\line K_{n_1,n_2})$.  Even in the
special case when $n_1=n_2=n$ (so $K(\line K_{n,n})$ is $n$-regular),
the exact sequence
$$
0 \rightarrow \underbrace{\ZZ_n^{\beta-2} \oplus
  \ZZ_n}_{\ZZ_n^{n(n-2)}} \longrightarrow \underbrace{K(\line K_{n,n})}_{\ZZ_{2n^2}^{2(n-2)} \oplus \ZZ_{2n}^{(n-2)^2+1}}
\overset{f}{\longrightarrow} \underbrace{K(\sd K_{n,n})}_{\ZZ_{2n^2}^1 \oplus \ZZ_{2n}^{2(n-2)} \oplus \ZZ_2^{(n-2)^2}} \longrightarrow
\ZZ_n \rightarrow 0
$$
from Theorem~\ref{regular-result} does not determine {\it a priori} Berget's formula for $K(\line K_{n,n})$.



However, we note that at least Theorem~\ref{divisibility-theorem} does
predict that the expression
$$
K(\line K_{n_1,n_2})=\bigoplus_{i=1}^{\beta} \ZZ_{e_i} \text{ where }
\beta:=\beta(K_{n_1,n_2}) = (n_1-1)(n_2-1)
$$
should have $|V|=n_1+n_2$ dividing every one of the factors $e_i$.
This follows from equation~\eqref{convenient-semiregular-info}, since
$K_{n_1,n_2}$ is bipartite $(n_1,n_2)$-semiregular.  Hence for each
prime $p$, one has
$$
K/p^{k(p)} K \cong \ZZ_{p^{k(p)}}^{\beta-1} \oplus \ZZ_{\gcd(p^{k(p)},|V|)} = 
\ZZ_{p^{k(p)}}^{\beta},
$$
where $k(p)$ is the largest power such that $p^{k(p)}$ divides $n_1+n_2$.
Hence $K/(n_1+n_2)K \cong \ZZ_{n_1+n_2}^\beta$.

\subsection{The $d$-dimensional cube}
\label{cube-section}

Denote by $G_{\dcube}$ the graph of vertices and edges in the
$d$-dimensional cube, that is, $G_{\dcube}=(V,E)$ in which $V$ is the
set of all binary strings of length $d$, and $E$ has an edge between
any two such strings that differ in exactly one binary digit. This is
a $d$-regular bipartite graph, having
$$
\beta:=\beta(G_{\dcube})=(d-2)2^{d-1}+1.
$$  
One knows its spanning tree number (see, e.g., \cite[Example 5.6.10]{Stanley}):
$$
  \kappa(G_{\dcube}) = \frac{1}{2^d} \prod_{k=1}^d (2k)^{\binom{d}{k}} = 2^{2^d-d-1} \prod_{k=2}^d k^{\binom{d}{k}}.
$$
Correspondingly, work of H. Bai \cite{Bai} computes its critical group
structure away from the prime $2$: For odd primes $p$, one has
$$
\Syl_p (K(G_{\dcube})) = \Syl_p \left( \bigoplus_{k=2}^d
  \ZZ_k^{\binom{d}{k}} \right).
$$
Unfortunately, $\Syl_2 (K(G_{\dcube}))$ is a $2$-group that is still
not known for all $d$.

Consequently, Proposition~\ref{Lorenzini's-theorem} shows that
$K(\sd(G_{\dcube}))$ has the same $p$-primary structure as
$K(G_{\dcube})$ for odd primes $p$, and Theorem~\ref{regular-result}
gives the following exact sequence for every odd prime $p$:
\begin{equation}
\label{cube-at-odd-prime-sequence}
0 \rightarrow 
\Syl_p(\ZZ_d^{\beta-1}) \rightarrow
\Syl_p (K(\line G_{\dcube})) \rightarrow
\Syl_p \left( \bigoplus_{k=2}^d \ZZ_k^{\binom{d}{k}} \right) \rightarrow 
\ZZ_d \rightarrow 0.
\end{equation}
This is particularly effective when $d$ itself is prime since then
$\Syl_d \left( \bigoplus_{k=2}^d \ZZ_k^{\binom{d}{k}} \right) =
\ZZ_d$ and the exact sequence \eqref{cube-at-odd-prime-sequence} implies
that for odd primes $p$,
$$
\Syl_p (K(\line G_{\dcube})) = 
\begin{cases}
\ZZ_d^{\beta-1} &\text{ if }p=d, \\
\Syl_p \left( \bigoplus_{k=2}^{d-1} \ZZ_k^{\binom{d}{k}} \right) &\text{ if }p \neq d.
\end{cases}
$$
Meanwhile $\Syl_2 (K(\line G_{\dcube})) = \Syl_2 (K(\sd(G_{\dcube})))$
is unknown, but by Proposition~\ref{Lorenzini's-theorem}, is
completely determined by the unknown $2$-group
$\Syl_2(K(G_{\dcube}))$.

\subsection{The Platonic solids}
\label{Platonic-solid-section}

One source of regular graphs are the $1$-skeleta (= graph of vertices
and edges) of the Platonic solids.  There are certain features that
apply to any graph $G_P$ which is the $1$-skeleton of a
$3$-dimensional polyhedron $P$, and hence to any Platonic solid:
\begin{enumerate}  
\item[$\bullet$] Because the cycles surrounding the (polygonal) faces
  of $P$ generate the cycle lattice $Z$, the graph $G_P$ is bipartite
  if and only if each face of $P$ is an even $n$-gon.
\item[$\bullet$] Furthermore, the cycles that bound all but one face
  of $P$ form a basis for $Z$, so that $\beta(G_P)$ is always one less
  than the number of faces.
\item[$\bullet$] Such graphs $G_P$ are always $2$-edge-connected, so
  that Theorem~\ref{2-edge-connected-theorem} always applies.
\item[$\bullet$] Dual polyhedra $P, P^*$ have $G_P, G_{P^*}$ dual as
  planar graphs.  This identifies the lattice of bonds for one with
  the lattice of cycles for the other, and implies that their critical
  groups $K(G_P), K(G_{P*})$ are isomorphic; see also
  \cite{CoriRossin}.
\end{enumerate}

\subsubsection{The tetrahedron}
The tetrahedron has $1$-skeleton $G_{\tetrahedron}=K_4$, and hence
implicitly was discussed already in
Section~\ref{complete-graph-section} on $K_n$, as the special case
$n=4$.

\subsubsection{The cube and octahedron}
The cube and the octahedron are dual polyhedra.  Either by direct
computer calculation, or by noting $G_{\octahedron} \cong \line K_4$
and applying Corollary~\ref{complete-graph-corollary} with $n=4$, one
finds that
$$
\begin{aligned}
K(G_{\cube})=K(G_{\octahedron})&=\ZZ_2 \oplus \ZZ_8 \oplus \ZZ_{24} \\
                               &=\ZZ_2 \oplus \ZZ_8^2 \oplus \ZZ_3.
\end{aligned}
$$

Since $G_{\octahedron}$ is $4$-regular and nonbipartite with
$\beta(G_{\octahedron})=7$,
Corollary~\ref{regular-nonbipartite-corollary} then implies
$$
K(\line G_{\octahedron})
= \ZZ_2^2 \oplus \ZZ_8^2 \oplus \ZZ_{16} \oplus \ZZ_{64} \oplus \ZZ_{192}.
$$

For $G_{\cube}$, which has $\beta(G_{\cube})=5$, the results of
Section~\ref{cube-section} apply, and are particularly effective
because $d=3$ is prime.  They show that $ \Syl_p (K(\line G_{\dcube}))
$ vanishes except for $p=2,3$, with
\[
\begin{aligned}
\Syl_3 (K(\line G_{\dcube})) &= \ZZ_3^4\\
\Syl_2 (K(\line G_{\dcube})) 
      &= \Syl_2 (K(\sd(G_{\dcube}))) 
       = \ZZ_2^2 \oplus \ZZ_4 \oplus \ZZ_{16}^2.
\end{aligned}
\]
Hence 
\[
\begin{aligned}
K(\line G_{\cube} ) &= \ZZ_2^2 \oplus \ZZ_4 \oplus \ZZ_{16}^2 \oplus \ZZ_3^4\\
                    &= \ZZ_2 \oplus \ZZ_6 \oplus \ZZ_{12} \oplus \ZZ_{48}^2.
\end{aligned}
\]

\subsubsection{The dodecahedron and icosahedron}
The dodecahedron and icosahedron are dual polyhedra, both of whose
graphs are nonbipartite.  Computer calculation shows that
$$
\begin{aligned}
K(G_{\dodecahedron})=K(G_{\icosahedron}) &= \ZZ_2 \oplus \ZZ_{12} \oplus \ZZ_{60}^3  \\
                                         &= \ZZ_2 \oplus \ZZ_4^4 \oplus \ZZ_3^4 \oplus \ZZ_5^3.
\end{aligned}
$$
Since $G_{\dodecahedron}$ is $5$-regular with $\beta(G_{\dodecahedron})=11$, 
one concludes from
Corollary~\ref{regular-nonbipartite-corollary} that
$$
\begin{aligned}
K(\line G_{\dodecahedron})
&= \ZZ_2^2 \oplus \ZZ_6^4 \oplus \ZZ_{12} \oplus \ZZ_{72} \oplus \ZZ_{360}^3\\
&=  \ZZ_2^6 \oplus \ZZ_4 \oplus \ZZ_8^4 \oplus \ZZ_3^5 \oplus \ZZ_9^4 \oplus \ZZ_5^3. \\
\end{aligned}
$$
Since $G_{\icosahedron}$ is $3$-regular with $\beta(G_{\icosahedron})=19$, 
one concludes from
Corollary~\ref{regular-nonbipartite-corollary} that
$$
\begin{aligned}
K(\line G_{\icosahedron}) 
& = \ZZ_2^2 \oplus \ZZ_{10}^{12} \oplus \ZZ_{20} \oplus \ZZ_{120} \oplus \ZZ_{600}^3\\
& = \ZZ_2^{14} \oplus \ZZ_4 \oplus \ZZ_8^4 \oplus \ZZ_3^4 \oplus \ZZ_5^{14} \oplus \ZZ_{25}^3. \\
\end{aligned}
$$

\section*{Acknowledgments}
The authors thank the anonymous referee for a careful reading of the
manuscript and many helpful suggestions, and also thank David Treumann
for allowing them to include some of his results here.


\end{document}